\documentclass{article}

\title{The Topological Fundamental Group and Free Topological Groups}
\author{Jeremy Brazas}
\usepackage{stmaryrd,graphicx}
\usepackage{amssymb,mathrsfs,amsmath,amscd,amsthm,fullpage}
\usepackage[all,cmtip]{xy}
\DeclareMathAlphabet{\mathpzc}{OT1}{pzc}{m}{it}
\usepackage{amsfonts,txfonts,pxfonts,latexsym,wasysym}
\newcommand{\sus}{\Sigma(X_{+})}
\newcommand{\ra}{\rightarrow}
\newcommand{\piz}{\pi_{0}(X)}

\newcommand{\piztop}{\pi_{0}^{top}(X)}
\newcommand{\pitopx}{\pi_{1}^{top}(X,x)}
\newcommand{\pitopxx}{\pi_{1}^{top}(X)}

\newcommand{\lpspx}{\Omega(X,x)}
\newcommand{\lpspxx}{\Omega(X)}
\newcommand{\slpsp}{\Omega_{+s}(\sus)}

\newcommand{\pitop}{\pi_{1}^{top}(\sus)}
\newcommand{\lpsp}{\Omega(\sus)}

\newcommand{\my}{M^{\ast}(Y)}
\newcommand{\jx}{M_{T}^{\ast}(X)}
\newcommand{\jy}{M_{T}^{\ast}(Y)}
\newcommand{\jtop}{M_{T}^{\ast}(\piztop)}
\newcommand{\jlpsp}{M_{T}^{\ast}(\Omega_{+s}(\sus))}
\newcommand{\mtpx}{M_{T}^{\ast}(P_X)}
\newcommand{\fry}{F_{R}(Y)}
\newcommand{\frx}{F_{R}(X)}
\newcommand{\ftx}{F_{M}(X)}
\newcommand{\fty}{F_{M}(Y)}
\newcommand{\frpi}{F_{R}\left(\piztop\right)}
\newcommand{\ftpi}{F_{M}\left(\piztop\right)}
\newcommand{\fotop}{F_{R}^{P_X}(\piztop)}
\newtheorem{theorem}{Theorem}[section]
\newtheorem{lemma}[theorem]{Lemma}
\newtheorem{proposition}[theorem]{Proposition}
\newtheorem{corollary}[theorem]{Corollary}
\newtheorem{definition}[theorem]{Definition}

\newtheorem{claim}[theorem]{Claim}
\newtheorem{example}[theorem]{Example}
\newtheorem{remark}[theorem]{Remark}

\begin{document}
\maketitle
\begin{abstract}
The topological fundamental group $\pi_{1}^{top}$ is a homotopy invariant finer than the usual fundamental group. It assigns to each space a quasitopological group and is discrete on spaces which admit universal covers. For an arbitrary space $X$, we compute the topological fundamental group of the suspension space $\Sigma(X_+)$ and find that $\pi_{1}^{top}(\Sigma(X_+))$ either fails to be a topological group or is the free topological group on the path component space of $X$. Using this computation, we provide an abundance of counterexamples to the assertion that all topological fundamental groups are topological groups. A relation to free topological groups allows us to reduce the problem of characterizing Hausdorff spaces $X$ for which $\pi_{1}^{top}(\Sigma(X_+))$ is a Hausdorff topological group to some well known classification problems in topology.
\end{abstract}

The fact that classical homotopy theory is insufficient for studying spaces with homotopy type other than that of a CW-complex has motivated the introduction of a number of invariants useful for studying spaces with complex local structure. For instance, in \u{C}ech theory, one typically approximates complicated spaces with "nice" spaces and takes the limit or colimit of an algebraic invariant evaluated on the approximating spaces. Another approach is to directly transfer topological data to algebraic invariants such as homotopy or (co)homology groups by endowing them with natural topologies that behave nicely with respect to the algebraic structure. While this second approach does not yield purely algebraic objects, it does have the advantage of allowing direct application of the rich theory of topological algebra. The notion of "topologized" homotopy invariant seems to have been introduced by Hurewicz in \cite{Hurewicz} and studied subsequently by Dugundji in \cite{D}. Whereas these early methods focused on "finite step homotopies" through open covers of spaces, we are primarily interested in the properties of a topologized version of the usual fundamental group.\\

The topological fundamental group $\pitopx$ of a based space $(X,x)$, as first sepecified by Biss \cite{Biss}, is the fundamental group $\pi_{1}(X,x)$ endowed with the natural topology that arrises from viewing it as a quotient space of the space of loops based at $x$. This choice of topological structure, makes $\pi_{1}^{top}$ particularly useful for studying the homotopy of spaces that lack universal covers, i.e. that fail to be locally path connected or semi-locally simply connected. Previously, several authors asserted that topological fundamental groups are always topological groups \cite{Biss,F4,GHMM,SS}, overlooking the fact that products of quotient maps are not always quotient maps. The initial intention of the research presented here was to produce couterexamples to this assertion. Recently, Fabel \cite{F3} has shown that the Hawaiian earring group $\pi_{1}^{top}(\mathbb{HE},(0,0))$ fails to be a topological group and here we find an abundance of spaces whose topological fundamental group fails to be a topological group.\footnote[1]{
One may similarly define the topological fundamental groupoid \cite{M1} of an unbased space, however, the same care must be taken with respect to products of quotient maps.} The production of such spaces, however, has led to surprising connections to the study of free topological groups. Surprisingly, multiplication can fail to be continuous even for a space as nice as a locally simply connected planar continuum (see Example \ref{rationalexample}).\\

The existence of such examples then begs the question: What type of object is $\pitopx$? In section 1, we find that topological fundamental groups are quasitopological groups in the sense of \cite{AT}. Here, we also provide preliminaries and some of the basic theory of topological fundamental groups.\\

The counterexamples mentioned above come from a class of spaces considered in section 2. This class consists of reduced suspensions of spaces with disjoint basepoint (written $\sus$). The well known suspension-loop adjunction then provides an unexpected relation to the free (Markov) topological groups. In fact, $\pitop$ either fails to be a topological group or is the free topological group on the path component space of $X$. This new connection is particularly surprising since, in general, it is difficult to describe the topological structure of both topological fundamental groups and free topological groups. We also note that much has been done to determine when free topological groups have certain quotient structures and the realization of these objects as homotopy invariants may indicate a potential application to their study. The second purpose of this paper, and our main result, is to fully describe the isomorphism class of $\pitop$ in the category of quasitopological groups. Specifically, $\pitop$ is the quotient (via reduction of words) of the path component space of the free topological monoid on $X\sqcup X^{-1}$.\\

Determining when a topology on a group is a group topology (is such that multiplication and inversion are continuous) is fundamental to the theory of topological groups. In section 3, we apply this theory to describe the topological properties of $\pitop$ and reduce the classification of Hausdorff spaces $X$ such that $\pitop$ is a Hausdorff topological group (and necessarily a free topological group) to three separate and well known classification problems in topology. We find that $\pitop$ is a Hausdorff topological group if and only if all four of the following conditions hold:
\begin{enumerate}
\item The path component space $\piztop$ is Hausdorff.
\item All powers of the path component quotient map $P_X:X\ra \piztop$ are again quotient maps.
\item The free topological group on the path component space $\piztop$ has the inductive limit topology of subspaces $F_{n}(\piztop)$ consisting of words of length $\leq n$.
\item For every $n\geq 1$, the canonical multiplication map $\mathbf{i}_{n}:\coprod_{i=0}^{n}(\piztop\sqcup \piztop^{-1})^{i}\ra F_{n}(\piztop)$ is a quotient map.
\end{enumerate}
Determining when products of quotient maps are quotient maps is fundamental to topology and has been extensively studied. Additionally, the classification of spaces for which conditions 3. and 4. are individually held are important open problems in the study of free topological groups.
\section{The Topological Fundamental Group}
We now recall some preliminaries regarding path component spaces and the compact-open topology.
\subsection{Path component spaces}
The path component space of a topological space $X$ is the set of path components $\pi_{0}(X)$ of $X$ with the quotient topology with respect to the canonical map $P_{X}:X\ra \pi_{0}(X)$. We denote this space as $\piztop$ and remove or change the subscript of the map $P_X$ when convenient. Since continuous functions take path components into path components,
a map $f:X\ra Y$ induces a map $f_{\ast}:\piztop\ra \pi_{0}^{top}(Y)$ taking the path component of $x$ in $X$ to the path component of $f(x)$ in $Y$. To see that $f_{\ast}$ is continuous, we observe that $P_{X}^{-1}(f_{\ast}^{-1}(U))=f^{-1}(P_{Y}^{-1}(U))$ is open in $X$ for all open $U\subseteq \pi_{0}^{top}(Y)$. Since $P_X$ is a quotient map, $f_{\ast}^{-1}(U)$ is open in $\piztop$.\\

If $X$ has basepoint $x$, we choose the basepoint of $\piztop$ to be the path component of $x$ in $X$. This gives an unbased and based version of the functor $\pi_{0}^{top}$, however, the presence of basepoint will be clear from context. The following remarks illustrate some of the properties of the path component space functor.
\begin{definition}
A space $X$ is \textbf{0-semilocally simply connected (0-SLSC)} if for each point $x\in X$, there is an open neighborhood $U$ of $x$ such that the inclusion $i:U\hookrightarrow X$ induces the constant map $i_{\ast}:\pi_{0}^{top}(U)\ra \piztop$.
\end{definition}
\begin{remark} \label{piztopdiscrete}
$X$ is 0-SLSC if and only if $\piztop$ has the discrete topology.
\end{remark}
\begin{remark} \label{quotientpathcomponent}
$\pi_{0}^{top}$ preserves coproducts and quotients, but does not preserve products. The non-topological path component functor $\pi_{0}:Top\ra Set$, of course, does preserve products.
\end{remark}
\begin{remark} \label{essentialsurj}
Every topological space $Y$ is homeomorphic to the path component space of some paracompact Hausdorff space $\mathcal{H}(Y)$. Some properties and variants of the functor $\mathcal{H}$ are included in \cite{Harris}.
\end{remark}
\subsection{The compact-open topology}
For spaces $X,Y$, let $M(X,Y)$ be the space of unbased maps $X\ra Y$ with the compact-open topology. A subbasis for this topology consists of neighborhoods of the form $\langle C,U\rangle=\{f|f(C)\subset U\}$ where $C\subseteq X$ is compact and $U$ is open in $Y$. If $A\subseteq X$ and $B\subseteq Y$, then $M(X,A;Y,B)\subseteq M(X,Y)$ is the subspace of relative maps (ie. such that $f(A)\subseteq B$) and if $X$ and $Y$ have basepoints $x$ and $y$ respectively, $M_{\ast}(X,x;Y,y)$ (or just $M_{\ast}(X,Y)$) is the subspace of basepoint preserving maps. In particular, $\Omega(X,x)=M_{\ast}(S^1,(1,0);X,x)$ is the space of based loops. When $X$ is path connected and the basepoint is clear, we just write $\Omega(X)$. For convenience, we will often replace $\Omega(X,x)$ by the homeomorphic relative mapping space $M(I,\{0,1\};X,\{x\})$ where $I=[0,1]$ is the closed unit interval. We will consistently denote the constant path at $x$ by $c_x:I\ra X$.\\

For any fixed, closed subinterval $A\subseteq I$, we make use of the following notation. Let $H_A:I\ra A$ be the unique, increasing, linear homeomorphism. For a path $p:I\ra X$, the restricted path of $p$ to $A$ is the composite $p_{A}=p|_{A}\circ H_{A}:I\ra A\ra X$. For integer $n\geq 1$ and $j=1,...,n$, let $K_{n}^{j}$ be the closed subinterval $\left[\frac{j-1}{n},\frac{j}{n}\right]\subseteq I$. If $p_j:I\ra X$, $j=1,...,n$ are paths such that $p_j(1)=p_{j+1}(0)$ for each $j=1,...,n-1$, then the n-fold concatenation of these paths is the unique path $q=\ast_{j=1}^{n}p_i=p_1\ast p_2\ast \dots\ast p_n$ such that $q_{K_{n}^{j}}=p_j$ for each $j$. It is well known that concatenation $\ast:\lpspxx\times \lpspxx\ra\lpspxx$ and loop inversion $^{-1}:\lpspxx\ra \lpspxx$, $p^{-1}(t)=p(1-t)$ are continuous with respect to the compact-open topology.\\

If $\mathscr{U}=\bigcap_{j=1}^{n}\langle C_j,U_j\rangle$ is a basic open neighborhood of a path $p\in M(I,X)$, then $\mathscr{U}_{A}=\bigcap_{A\cap C_j\neq \emptyset}\langle H_{A}^{-1}(A\cap C_j),U_j\rangle$ is a basic open neighborhood of $p_{A}$ called the \textbf{restricted neighborhood} of $\mathscr{U}$ to $A$. On the other hand, if $\mathscr{U}=\bigcap_{j=1}^{n}\langle C_j,U_j\rangle$ is a basic open neighborhood of the restricted path $p_{A}$, then $\mathscr{U}^{A}=\bigcap_{j=1}^{n}\langle H_{A}(C_j),U_j\rangle$ is a basic open neighborhood of $p$ called the \textbf{induced neighborhood} of $\mathscr{U}$ on $A$. For all basic open neighborhoods $\mathscr{U}$ in $M(I,X)$ we have that $$(\mathscr{U}^{A})_{A}=\mathscr{U}\text{  and  }\mathscr{U}\subseteq (\mathscr{U}_{A})^{A}$$ and if $t\in (0,1)$, then $$\mathscr{U}=(\mathscr{U}_{[0,t]})^{[0,t]}\cap (\mathscr{U}_{[t,1]})^{[t,1]}.$$

The following lemma is a basic property of free path spaces. For details see \cite{BR}.
\begin{lemma} \label{basis}
Suppose $\mathscr{B}_{X}$ is a basis for the topology of $X$ which is closed under finite intersection. The collection of open neighborhoods of the form $\bigcap_{j=1}^{n}\langle K_{n}^{j},U_j\rangle$ where $U_j\in \mathscr{B}_{X}$ is a basis for the compact-open topology of the free path space $M(I,X)$. Moreover, this basis is closed under finite intersection. 
\end{lemma}

This allows us to intuit basic open neighborhoods of paths and loops as finite, ordered sets of "instructions."
\subsection{The topological fundamental group and functorality}

The topological fundamental group of a based space $(X,x)$ is the path component space $\pitopx=\pi_{0}^{top}(\lpspx)$. For an unbased space $Y$, let $Y_+=Y\sqcup \{\ast\}$ denote the based space with added disjoint basepoint. The natural homeomorphisms $M_{\ast}(\Sigma(I_+),X)\cong M(I,\Omega(X))$ indicate that homotopy classes of based loops in $X$ are the same as path components in $\lpspx$. Consequently, $\pitopx$ may be described as the usual fundamental group with the final topology with respect to the canonical map $\lpspx\ra \pi_{1}(X,x)$ identifying homotopy classes of loops. 
\begin{proposition} \label{finesttopology}
The topology of $\pitopx$ is the finest topology on $\pi_{1}(X,x)$ such that the canonical map $\lpspx\ra \pi_{1}(X,x)$ is continuous.
\end{proposition}

While $\pitopx$ is a group with a topology, we show in section 3 that it is not always a topological group. With this in mind, we describe the functorial nature of $\pi_{1}^{top}$ using the next definition.
\begin{definition} \emph{
A \textbf{quasitopological group} $G$ is a group with topology such that inversion $G\ra G$, $g\mapsto g^{-1}$ is continuous and multiplication $G\times G\ra G$ is continuous in each variable (i.e. all translations are continuous). A morphism of quasitopological groups is a continuous homomorphism and the category of quasitopological groups is denoted $\mathscr{QTG}$.
}
\end{definition}

A basic account of the theory of quasitopological groups may be found in \cite{AT}. The next two lemmas combine results from \cite{Biss} and \cite{CM}.
\begin{lemma} \label{pi1topfunctorality} If $hTop_{\ast}$ denotes the homotopy category of based topological spaces, then $\pi_{1}^{top}:hTop_{\ast}\ra \mathscr{QTG}$ is a functor.
\end{lemma}
\begin{proof} It was noted in section 1.2 that both loop inversion $\lpspx\ra \lpspx$, $\alpha\mapsto \alpha^{-1}$ and concatenation $\lpspx\times \lpspx\ra \lpspx$, $(\alpha,\beta)\mapsto \alpha\ast\beta$ are continuous. For any fixed loop $\alpha\in \lpspx$, concatenation restricts to the continuous maps $r_{\alpha},l_{\alpha}:\lpspx\ra \lpspx$ which are right and left concatenation by $\alpha$ respectively. Applying the path component space functor to loop inversion, $r_{\alpha}$, and $l_{\alpha}$, it is clear that the functions $\pitopx\ra \pitopx$ given by group inversion, and right and left multiplication by $[\alpha]$ are continuous. Therefore, $\pitopx$ is a quasitopological group. For each based map $f:(X,x)\ra (Y,y)$, there is an induced homomorphism $f_{\ast}:\pi_{1}(X,x)\ra \pi_{1}(Y,y)$ on fundamental groups. By the functorality of $\pi_{0}^{top}$ and $\Omega$, $f_{\ast}=\pi_{0}^{top}(\Omega(f))$ is also continuous. Preservation of identity and composition are given by functorality of $\pi_{1}$. Since two homotopic maps induce the same continuous homomorphism, it follows that $\pi_{1}^{top}$ is well defined on the homotopy category.
\end{proof}
\begin{lemma} \label{basepoint}
If $x$ and $y$ lie in the same path component of $X$, then $\pitopx\cong \pi_{1}^{top}(X,y)$ as quasitopological groups.
\end{lemma}
\begin{proof}
Suppose $\gamma:I\ra X$ is a path with $\gamma(0)=x$ and $\gamma(1)=y$. The maps $\Gamma:\Omega(X,y)\ra \Omega (X,x)$, $\Gamma(\alpha)=\gamma\ast\alpha\ast\gamma^{-1}$ and $\Gamma ':\Omega(X,x)\ra \Omega (X,y)$, $\Gamma ' (\beta)=\gamma^{-1}\ast\beta\ast\gamma$ are continuous and induce the appropriate continuous inverses on path component spaces.
\end{proof}

Since choice of basepoint within path components is irrelivant, we suppress the basepoint and simply write $\pitopxx$ whenever $X$ is path connected.\\

While multiplication is not always continuous in topological fundamental groups, we do have the continuity of power maps.
\begin{proposition}
For each integer $n\geq 1$, the n-th power map $pow_n:\pitopx\ra \pitopx$, $pow_n([\alpha])=[\alpha]^n$ is continuous.
\end{proposition}
\begin{proof}
Let $\Delta_n:\lpspx\ra \lpspx^{n}$, $\Delta_n(\alpha)=(\alpha,...,\alpha)$ be the diagonal map and $m_n:\lpspx^{n}\ra \lpspx$, $m_n(\alpha_1,...,\alpha_n)=\alpha_1\ast\dots \ast\alpha_{n}$ be the n-fold concatenation map. The composite $m_n\circ \Delta_n:\lpspx\ra \lpspx$ is continuous and when we apply the path component space functor, we get a continuous map $$\pitopx\ra \pitopx \text{  where  }[\alpha]\ra \left[\ast_{i=1}^{n}\alpha\right]=[\alpha]^{n}$$which is precisely the n-th power map.
\end{proof}
It is well known that there are quasitopological groups with discontinuous power maps. For instance, consider subset $K=\left\{\frac{\epsilon}{3^n}|n\geq 1,\epsilon=\pm 1\right\}$ of the additive group of reals $\mathbb{R}$. Form a neighborhood base at $0$ consisting of symmetric neighborhoods of the form $(-t,t)-K$. Letting all translations be homemeorphisms we generate a topology which makes $\mathbb{R}$ a quasitopological group. It is not too difficult to see that the square map $s:\mathbb{R}\ra \mathbb{R}$, $s(t)=2t$ is discontinuous. In particular, the sequence $\frac{1}{2(3^n)}$ converges to $0$ but $s\left(\frac{1}{2(3^n)}\right)=\frac{1}{3^n}$ does not. The existence of such groups illustrates another defect of $\pi_{1}^{top}$, namely that $\pi_{1}^{top}$ is not essentially surjective:
\begin{corollary}
Not every quasitopological group is a topological fundamental group.
\end{corollary}
\subsection{A characterization of discreteness}

In general, it is difficult to determine if the topological fundamental group of a space is a topological group. There are, however, instances when it is easy to answer in the affirmative, namely those spaces $X$ for which $\pitopxx$ has the discrete topology.

\begin{lemma}
For any based space $(X,x)$, the following are equivalent:
\begin{enumerate}
\item $\pitopx$ has the discrete topology.
\item $\lpspx$ is 0-SLSC.
\item The singleton $\{[c_x]\}$ containing the identity is open in $\pitopx$.
\item Each null-homotopic loop $\alpha\in\lpspx$ has an open neighborhood containing only null-homotopic loops.
\end{enumerate}
\end{lemma}
\begin{proof} 1. $\Leftrightarrow$ 2. follows from Remark \ref{piztopdiscrete}. 1. $\Leftrightarrow$ 3. follows from the fact that left and right translations in $\pitopx$ are homeomorphisms. 3. $\Leftrightarrow$ 4. is due to the fact that $\pitopxx$ has the quotient topology of $\lpspx$.\end{proof}

These obvious characterizations are inconvenient in that they do not characterize discreteness in terms of the topological properties of $X$ itself. We give such a description in the next theorem which has also appeared in \cite{CM}. Recall that a space is \textbf{semilocally simply connected (SLSC)} if for each $x\in X$ there is an open neighborhood $U$ of $x$ such that the inclusion $i:U\hookrightarrow X$ induces the trivial homomorphism $i_{\ast}:\pi_{1}(U,x)\ra \pi_{1}(X,x)$. Note that if $i:U\hookrightarrow X$ induces the trivial homomorphism $i_{\ast}:\pi_{1}(U,x)\ra \pi_{1}(X,x)$ and $z$ lies in the same path component of $U$ as $x$, then the induced homomorphism $i_{\ast}:\pi_{1}(U,z)\ra \pi_{1}(X,z)$ is also trivial.
\begin{theorem} \label{discreteness}
Suppose $X$ is path connected. If $\pitopxx$ is discrete, then $X$ is SLSC. If $X$ is locally path connected and SLSC, then $\pitopxx$ is discrete.
\end{theorem}
\begin{proof} We suppose $x\in X$ and by Lemma \ref{basepoint} may assume that $\pitopx$ is discrete or equivalently that $\lpspx$ is 0-SLSC. This allows us to find open neighborhood $W$ of the constant loop $c_x$ in $\lpspx$ such that $\alpha\simeq c_x$ for each $\alpha\in W$. There is an open neighborhood $U$ of $x$ in $X$ such that $c_x\in \langle S^1,U\rangle\subseteq W$. Since every loop $\alpha\in\langle S^1,U\rangle$ is null-homotopic in $X$, the inclusion $i:U\hookrightarrow X$ induces the trivial homomorphism $i_{\ast}:\pi_{1}(U,x)\ra \pi_{1}(X,x)$. Thus $X$ is SLSC.\\
\\
We now suppose $X$ is locally path connected and SLSC and that $\alpha\in M(I,\{0,1\};X,\{x\})$. We find an open neighborhood of $\alpha$ in $M(I,\{0,1\};X,\{x\})$ containing only loops homotopic to $\alpha$ in $X$. This suffices to show that $\Omega(X,x)$ is 0-SLSC. For each $t\in I$, we find an open neighborhood $U_t$ of $\alpha(t)$ in $X$ such that the inclusion $u_t:U_t\hookrightarrow X$ induces the trivial homomorphism $(u_t)_{\ast}:\pi_{1}(U_t,\alpha(t))\ra \pi_{1}(X,\alpha(t))$. We then find a path connected, open neighborhood $V_t$ of $\alpha(t)$ contained in $U_t$. Take a finite subcover $\{V_{t_1},...,V_{t_k}\}$ of $\alpha(I)$ and finite subdivisions of $I$ to find an integer $m\geq 1$ such that $\alpha\in \bigcap_{j=1}^{m}\langle K_{m}^{j},V_j\rangle$ where $V_j=V_{t_{i_j}}$ for not necessarily distinct $i_j\in \{1,...,k\}$. For $j=0,...,m$, let $s_j=\frac{j}{m}\in I$. For each $j=1,...,m-1$ we have $\alpha(s_j)\in V_{j}\cap V_{j+1}$ and find a path connected, open neighborhood $W_j$ such that $\alpha(s_j)\in W_j\subseteq V_{j}\cap V_{j+1}$. Now$$\mathscr{U}=\bigcap_{j=1}^{m}\langle K_{m}^{j},V_j\rangle\cap \bigcap_{j=1}^{m-1}\langle \{s_j\}, W_j\rangle $$is an open neighborhood of $\alpha$ in $M(I,\{0,1\};X,\{x\})$. We suppose $\gamma\in \mathscr{U}$ and construct a homotopy to $\alpha$. We have $\gamma(s_j)\in W_j$ for $j=1,...,m-1$ allowing us to find paths $p_j:I\ra W_j$ such that $p_j(0)=\alpha(s_j)$ and $p_j(1)=\gamma(s_j)$. Let $p_0=p_m=c_x$ be the constant path at $x$. We now make use of our notation for restricted paths. For $j=1,...,m$ we define loops $\beta_j:I\ra V_j$ based at $\alpha(s_{j-1})$ as the concatenations $$\beta_j=p_{j-1}\ast\gamma_{K_{m}^{j}}\ast p_{j}^{-1}\ast \alpha_{K_{m}^{j}}^{-1}$$Recall that $V_j=V_{t_{i_j}}$ where $\alpha(t_{i_j})\in V_{j}$. Since $V_j$ is path connected, the points $\alpha(s_{j-1})$ and $\alpha(t_{i_j})$ lie in the same path component of $U_j$. Therefore the inclusion $u_j:U_j\hookrightarrow X$ induces the trivial homomorphism $(u_j)_{\ast}:\pi_{1}(U_j,\alpha(s_{j-1}))\ra \pi_{1}(X,\alpha(s_{j-1}))$. Consequently, each loop $\beta_j$ is homotopic (in $X$) to the constant loop at $\alpha(s_{j-1})$. The homotopies of loops $\beta_j\simeq c_{\alpha(s_{j-1})}$ give fixed endpoint homotopies of paths $\alpha_{K_{m}^{j}}\simeq p_{j-1}\ast\gamma_{K_{m}^{j}}\ast p_{j}^{-1}$. Now we have concatenations of homotopies$$\alpha\simeq \ast_{j=1}^{m}\alpha_{K_{m}^{j}}\simeq \ast_{j=1}^{m} \left(p_{j-1}\ast\gamma_{K_{m}^{j}}\ast p_{j}^{-1}\right)\simeq p_{0}\ast\left(\ast_{j=1}^{m}\gamma_{K_{m}^{j}}\right)\ast p_{m}^{-1}\simeq p_{0}\ast\gamma \ast p_{m}^{-1}\simeq \gamma$$This proves that $\mathscr{U}$ contains only loops homotopic to $\alpha$ in $X$, or in otherwords that the inclusion $\mathscr{U}\hookrightarrow M(I,\{0,1\};X,\{x\})$ induces the constant function on path components.\end{proof}
\subsection{Separation properties and a basis for $\pi_{1}^{top}(X)$}

In \cite{Biss} and \cite{F7}, the harmonic archepellago (a non-compact subspace of $\mathbb{R}^{3}$) is shown to have uncountable, indiscrete topological fundamental group. The next example is a compact metric space with topological fundamental group isomorphic to the indiscrete group of integers.
\begin{example} \emph{
Let $S^1=\{(x,y,0)\in\mathbb{R}^{3}|x^2+y^2=1\}$ be the unit circle in the xy-plane of $\mathbb{R}^{3}$. For all integers $n\geq 1$, we let $$C_n=\left\{(x,y,z)\in \mathbb{R}^{3}|\left(x-\frac{1}{n}\right)^{2}+y^2+z^2=\left(1+\frac{1}{n}\right)^{2}\right\}$$Now let $X=S^1\cup\left(\bigcup_{n\geq 1}C_n\right)$ have basepoint $(-1,0,0)$. This is a compact subspace of $\mathbb{R}^{3}$ and is weakly equivalent to the wedge of spheres $S^1\vee\left(\bigvee_{n\geq 1} S^{2}\right)$. We have $\pi_{1}(X)\cong \mathbb{Z}$, however, every open neighborhood $W$ of a loop $\alpha:S^1\ra S^1\subset X$ contains a loop $\beta:S^1\ra \bigcup_{n\geq 1}C_n\subset X$ which is null homotopic. Therefore every open neighborhood of the class $[\alpha]$ in $\pitopxx$ contains the identity. From this, one may easily use the properties of quasitopological groups to show that the topology of $\pitopxx$ is the indiscrete topology. This example also illustrates how weakly equivalent spaces may have non-isomorphic topological fundamental groups.}
\end{example}

Since there are simple spaces with non-trivial, indiscrete topological fundamental group, we cannot take any separation properties for granted. We now give a simple characterization of spaces with T1 topological fundamental group.
\begin{proposition} \label{T1lemma} Suppose $X$ is path connected. The following are equivalent:
\begin{enumerate}
\item For each loop $\alpha\in\Omega(X)$ which is not null-homotopic, there is an open neighborhood $V$ of $\alpha$ such that $V$ contains no null-homotopic loops.
\item The singleton containing the identity is closed in $\pitopxx$.
\item $\pitopxx$ is $T_0$.
\item $\pitopxx$ is $T_1$.
\end{enumerate}
\end{proposition}
\begin{proof} 1. $ \Leftrightarrow $ 2. follows from the definition of the quotient topology and 2. $\Leftrightarrow$ 3. $\Leftrightarrow $ 4. holds for all quasitopological groups.\end{proof}

As the continuity of multiplication is critical in proving that every $T_0$ topological group is Tychonoff, it can be difficult to recognize separation properties $T_i$, $i\geq 2$ in quasitopological groups. Additionally, the complex nature of homotopy as an equivalence relation further complicates our attempt to characterize stronger separation properties in topological fundamental groups. To be able to make any general statement for when $\pitopxx$ is Hausdorff, it will be necessary to find a basis of open neighborhoods at the identity of $\pitopxx$. While computationally challenging, there is a method of constructing a basis of open neighborhoods for any quotient space. We take this approach so that if $q:Y\shortrightarrow Z$ is a quotient map, a basis for $Z$ may be described in terms of open coverings of $Y$.
\begin{definition} \emph{ For any space $Y$, a \textbf{pointwise open cover} of $Y$ is an open cover $\mathscr{U}=\{U^{y}\}_{y\in Y}$ where each point $y\in Y$ has a distinguished open neighborhood $U^{y}$ containing it. Let $Cov(Y)$, be the directed set of pointed open covers of $Y$ where the direction is given by pointwise refinement: If $\mathscr{U}=\{U^{y}\}_{y\in Y},\mathscr{V}=\{V^{y}\}_{y\in Y}\in Cov(Y)$, then we say $\mathscr{U}\preceq \mathscr{V}$ when $V_y\subseteq U_y$ for each $y\in Y$.}
\end{definition}

We also make use of the following notation: If $\mathscr{U}=\{U^{y}\}_{y\in Y}\in Cov(Y)$ is a pointwise open covering of $Y$ and $A\subseteq Y$, let $\mathscr{U}(A)=\bigcup_{a\in A} U^{a}$.\\
\\
\textbf{Construction of basis:} Suppose $q:Y\shortrightarrow Z$ is quotient map, $z\in Z$, and $\mathscr{U}\in Cov(Y)$ is a fixed point-wise open cover. We construct open neighborhoods of $z$ in $Z$ in the most unabashed way, that is, by recursively "collecting" the elements of $Y$ so that our collection is both open and saturated. We begin by letting $\mathcal{O}_{q}^{0}\left(z,\mathscr{U}\right)=\{z\}$. For integer $n\geq 1$, we define $\mathcal{O}_{q}^{n}\left(z,\mathscr{U}\right)\subseteq Z$ as$$\mathcal{O}_{q}^{n}\left(z,\mathscr{U}\right)=q\left(\mathscr{U}\left(q^{-1}\left( \mathcal{O}_{q}^{n-1}\left(z,\mathscr{U}\right)\right)\right)\right)$$

It is clear that $\mathcal{O}_{q}^{n-1}\left(z,\mathscr{U}\right)\subseteq \mathcal{O}_{q}^{n}\left(z,\mathscr{U}\right)$ for all $n\geq 1$. We then may take the union $$\mathcal{O}_{q}\left(z,\mathscr{U}\right)=\bigcup_{n}\mathcal{O}_{q}^{n}\left(z,\mathscr{U}\right)$$

Note that if $y\in q^{-1}\left(\mathcal{O}_{q}\left(z,\mathscr{U}\right)\right)$, then $U^{y}\subset q^{-1}\left(\mathcal{O}_{q}\left(z,\mathscr{U}\right)\right)$ so that $\mathcal{O}_{q}\left(z,\mathscr{U}\right)$ is open in $Z$. Also, if $\mathscr{W}=\{W^{y}\}_{y\in Y}$ is another point-wise open cover of $Y$ such that $q(W^{y})\subseteq q(U^{y})$ for each $y\in Y$, then $\mathcal{O}_{q}\left(z,\mathscr{W}\right)\subseteq \mathcal{O}_{q}\left(z,\mathscr{U}\right)$. The neighborhood $\mathcal{O}_{q}\left(z,\mathscr{U}\right)$ is said to be the \textit{open neighborhood of} $z$ \textit{in} $Z$ \textit{generated by} $\mathscr{U}$. It is easy to see that for each open neighborhood $V$ of $z$ in $Z$, there is a pointwise open covering $\mathscr{V}\in Cov(Y)$ such that $z\in\mathcal{O}_{q}\left(z,\mathscr{V}\right)\subseteq V$. In particular, let $V^y=q^{-1}(V)$ when $q(y)\in V$ and $V^y=Y$ otherwise.\\

In the case of the quotient map $P:\Omega(X)\ra \pitopxx$, the neighborhoods $\mathcal{O}_{P}\left([c_x],\mathscr{U}\right)$ for $\mathscr{U}=\{U^\beta\}_{\beta\in\Omega(X)}\in Cov(\Omega(X))$ give a basis for the topology of $\pitopxx$ at the identity $[c_x]$. The loops in $P^{-1}\left(\mathcal{O}_{P}\left([c_x],\mathscr{U}\right)\right)$ can be described as follows:\\

For each $\alpha\in P^{-1}\left(\mathcal{O}_{P}\left([c_x],\mathscr{U}\right)\right)$, there is an integer $n\geq 1$ and a sequence of loops $\gamma_0,\gamma_1,...,\gamma_{2n+1}$ where $\gamma_0=c_x$, $\gamma_{2i}\simeq \gamma_{2i+1}$ for $i=0,1,...,n$, $\gamma_{2i+1}\in U^{\gamma_{2i+2}}$ for $i=0,1,...,n-1$, and $\gamma_{2n+1}=\alpha$. In this sense, the neighborhood $\mathcal{O}_{P}\left([c_x],\mathscr{U}\right)$ is an alternating "collection" of homotopy classes and nearby loops (the nearby being determined by the elements of $\mathscr{U}$).\\

We can state what it means for $\pitopxx$ to be Hausdorff in these terms.
\begin{proposition} $\pitopxx$ is Hausdorff if and only if for each class $[\beta]\in \pitopxx-\{[c_x]\}$, there is a pointwise open covering $\mathscr{U}\in Cov(\Omega(X))$ such that $\mathscr{O}_{P}([c_x],\mathscr{U})\cap \mathscr{O}_{P}([\beta],\mathscr{U})=\emptyset$.
\end{proposition}
\begin{proof} If $\pitopxx$ is Hausdorff and $[\beta]\in \pitopxx-\{[c_x]\}$, we can find disjoint open neighborhoods $W$ of $[c_x]$ and $V$ of $[\beta]$. Now we may find pointwise open coverings $\mathscr{W}=\{W^{\alpha}\}_{\alpha\in\Omega(X)},\mathscr{V}=\{V^{\alpha}\}_{\alpha\in\Omega(X)}\in Cov(\Omega(X))$ such that $\mathscr{O}_{P}([c_x],\mathscr{W})\subset W$ and $\mathscr{O}_{P}([\beta],\mathscr{V})\subset V$. We let $\mathscr{W}\cap\mathscr{V}=\{W^{\alpha}\cap V^{\alpha}\}_{\alpha\in\Omega(X)}\in Cov(\Omega(X))$ be the intersection of the two. Since $\mathscr{W},\mathscr{V}\preceq \mathscr{W}\cap\mathscr{V}$, we have $\mathscr{O}_{P}([c_x],\mathscr{W}\cap\mathscr{V})\subseteq \mathscr{O}_{P}([c_x],\mathscr{W})\subset W$ and $\mathscr{O}_{P}([\beta],\mathscr{W}\cap\mathscr{V})\subseteq \mathscr{O}_{P}([\beta],\mathscr{V})\subset V$.\\
\\
If the second statement holds, we suppose that $[\beta_1]$ and $[\beta_2]$ are distinct classes in $\pitopxx$. Therefore $[\beta_1\ast\beta_2^{-1}]\neq [c_x]$ and by assumption there is a $\mathscr{U}\in Cov(\Omega(X))$ such that $\mathscr{O}_{P}([c_x],\mathscr{U})\cap \mathscr{O}_{P}([\beta_1\ast\beta_2^{-1}],\mathscr{U})=\emptyset$. Since right multiplication by $[\beta_2]$ is a homeomorphism, we have that $\mathscr{O}_{P}([\beta_1\ast\beta_2^{-1}],\mathscr{U})[\beta_2]$ is open containing $[\beta_1]$ and $\mathscr{O}_{P}([c_x],\mathscr{U}) [\beta_2]$ is open containing $[\beta_2]$. But  $\left(\mathscr{O}_{P}([\beta_1\ast\beta_2^{-1}],\mathscr{U})[\beta_2]\right)\cap \left( \mathscr{O}_{P}([c_x],\mathscr{U}) [\beta_2]\right)=\emptyset$ and so $\pitopxx$ is Hausdroff.\end{proof}
\section{A Computation of $\pitop$}

We now describe the isomorphim class of $\pitop$ in the category of quasitopological groups for an arbitrary space $X$.
\subsection{The spaces $\sus$}

Suppose $X$ is a topological space and $X_+=X\sqcup\{\ast\}$ is the based space with added disjoint basepoint. Let $$(\sus,x_0)=\left(\frac{X_{+}\times I}{X\times \{0,1\}\cup \{\ast\}\times I},x_0\right)$$ be the reduced suspension of $X_+$ with canonical choice of basepoint and $x\wedge s$ denote the image of $(x,s)\in X\times I$ under the quotient map $X_{+}\times I\ra \sus$. For subsets $A\subseteq X$ and $S\subseteq I$, let $A\wedge S=\{a\wedge s|a\in A,s \in S\}$. A subspace $P\wedge I$ where $P\in\pi_{0}(X)$ is a path component of $X$ is called a \textbf{hoop} of $\sus$.\\

Suppose $\mathscr{B}_{X}$ is a basis for the topology of $X$ which is closed under finite intersections. For a point $x\wedge t\in X\wedge (0,1)=\sus-\{x_0\}$, a subset $U\wedge (c,d)$ where $x\in U$, $U\in \mathscr{B}_{X}$ and $t\in(c,d)\subseteq (0,1)$ is an open neighborhood of $x\wedge t$. Open neighborhoods of $x_0$ may be given in terms of open coverings of $ X\times \{0,1\}$ in $X\times I$. If $U^x\in\mathscr{B}_{X}$ is an open neighborhood of $x$ in $X$ and $t_x\in (0,\frac{1}{8})$, the set $$\bigcup_{x\in X} \left(U^x\wedge [0,t_x)\cup(1-t_x,1]\right)$$ is an open neighborhood of $x_0$ in $\sus$. The collection $\mathscr{B}_{\sus}$ of neighborhoods of the form $U\wedge (c,d)$ and $\bigcup_{x\in X}\left(U^x\wedge  [0,t_x)\cup(1-t_x,1]\right)$ is a basis for the topology $\sus$ which is closed under finite intersection.\\

We remark on some of the basic topological properties of $\sus$.
\begin{remark} 
\label{susfacts}
For an arbitrary space $X$,
\begin{enumerate}
\item $\sus$ is path-connected.
\item $\sus-\{x_0\}=X\wedge (0,1)\cong  X\times(0,1)$.
\item Every basic neighborhood $V\in \mathscr{B}_{\sus}$ containing $x_0$ is arc connected and simply connected.
\item For each $t\in (0,1)$, the closed subspaces $X\wedge [0, t]$ and $X\wedge [t,1]$ are homeomorphic to $CX$ the cone of $X$, and are contractible to the basepoint point $x_0$.
\item $\sus$ is Hausdorff if and only if $X$ is Hausdorff, but the following holds for arbitrary $X$: For each point $x\wedge t\in X \wedge (0,1)$, there are disjoint open neighborhoods separating $x\wedge t$ and the basepoint $x_0$.
\end{enumerate}
\end{remark}

The next remark illustrates that the spaces $\sus$ are natural generalizations of a wedge of circles. Intuitively, one might think of $\sus$ as a "wedge of circles parameterized by the space $X$."\\

Let $\bigvee_{X}S^{1}$ be the wedge of circles indexed by the underlying set of $X$. Suppose $\epsilon:I\ra S^1$ is the exponential map and a point in the $x$-th summand of the wedge is denoted as $\epsilon(t)_{x}$ for $t\in I$. The pushout property implies that every map $f:X\ra Y$ induces a map $\bigvee_{X}S^{1}\ra \bigvee_{Y}S^{1}$ given by $\epsilon(t)_x\mapsto \epsilon(t)_{f(x)}$ for all $t\in I$, $x\in X$. It is easy to see that $\bigvee_{(-)}S^{1}:Top\ra Top_{\ast}$ is a functor which we may relate to $\Sigma((-)_+)$ in the following way.
%
\begin{remark} \label{wedgegeneralizationlemma}
There is a natural transformation $\gamma:\bigvee_{(-)}S^{1}\ra \Sigma((-)_{+})$ where each component $\gamma_{X}:\bigvee_{X}S^{1}\ra \sus$ given by $\gamma_{X}(\epsilon(t)_x)=x\wedge t$ is a continuous bijection. Moreover, $\gamma_{X}$ is a homeomorphism if and only $X$ has the discrete topology.
\end{remark}

This fact gives some trivial cases for the topology of $\pitop$. Namely, if $X$ has the discrete topology, then $\sus$ is homeomorphic to a wedge of circles and by the Van-Kampen Theorem and Theorem \ref{discreteness}, $\pitop$ must be isomorphic to the discrete free group $F(X)$ on the underlying set of $X$. We will see later on that $\pitop$ is discrete if and only if $\piztop$ is discrete.
\begin{remark} \label{adjunction}
It is a well known fact that the reduced suspension functor $\Sigma:Top_{\ast}\ra Top_{\ast}$ is left adjoint to the loop space functor $\Omega:Top_{\ast}\ra Top_{\ast}$. Additionally, adding disjoint basepoint to an unbased space $(-)_+:Top\ra Top_{\ast}$ is left adjoint to the functor $U:Top_{\ast}\ra Top$ forgetting basepoint. Taking composites, we see the construction $\Sigma((-)_{+}):Top\ra Top_{\ast}$ is a functor. For a map $f:X\ra Y$, the map $\Sigma(f_{+}):\sus\ra \Sigma(Y_+)$ is defined by $\Sigma(f_{+})(x\wedge s)=f(x)\wedge s$. Moreover, $\Sigma((-)_{+})$ is left adjoint to the composite functor $U\Omega:Top_{\ast}\ra Top$. This is illustrated by the natural homeomorphisms $$M_{\ast}(\sus,Y)\cong M_{\ast}(X_+,\Omega(Y))\cong M(X,U\Omega(Y)).$$
\end{remark}

This adjunction immediately gives the following motivation for our proposed computation of $\pitop$. We will say that a quotient group $G/N$ of a quasitopological group $G$ is a topological quotient group when $G/N$ has the quotient topology with respect to the canonical homomorphism $G\ra G/N$.
\begin{proposition} \label{quotientofpitop}
Every topological fundamental group $\pi_{1}^{top}(Y)$ is the topological quotient of $\pitop$ for some space $X$.
\end{proposition}
\begin{proof} Let $cu:\Sigma(\Omega(Y)_{+})\ra Y$ be the adjoint of the unbased identity of $\Omega(Y)$. The basic property of counits gives that the unbased map $U\Omega(cu):\Omega(\Sigma(\Omega(Y)_{+}))\ra \Omega(Y)$ is a topological retraction. Applying the path component functor, we obtain a group epimorphism $\pi^{top}_{1}(\Sigma(\Omega(Y)_{+}))\ra \pi_{1}^{top}(Y)$ which is, by Remark \ref{quotientpathcomponent}, a quotient map of spaces. Take $X=\Omega(Y)$.\end{proof}
\subsection{Free topological monoids and the James map}

Prior to computing $\pi_{1}(\sus)$ and $\pitop$, we recall some common constructions in algebraic topology and topological algebra.\\

Let $Y^{-1}$ denote a homeomorphic copy of an unbased space $Y$. The free topological monoid on the disjoint union $Y\sqcup Y^{-1}$ is $$\jy=\coprod_{n\geq 0}(Y\sqcup Y^{-1})^{n}\text{ where }(Y\sqcup Y^{-1})^0=\{e\}.$$A non-empty element of $\jy$ is a word $y_{1}^{\epsilon_1}y_{2}^{\epsilon_2}\dots y_{n}^{\epsilon_n}$ where $\epsilon_i\in \{\pm 1\}$ and $y_{i}^{\epsilon_i}\in Y^{\epsilon_i}$. Multiplication is concatenation of words and the identity is the disjoint empty word $e$. The underlying monoid of $\jy$ is $\my$. The length of a word $w=y_{1}^{\epsilon_1}y_{2}^{\epsilon_2}\dots y_{n}^{\epsilon_n}$ is $|w|=n$ and we set $|e|=0$. For each finite (possibly empty) sequence $\zeta=\{\epsilon_i\}_{i=1}^{n}$ where $\epsilon_i\in \{\pm 1\}$, let $Y^{\zeta}=\{y_{1}^{\epsilon_1}y_{2}^{\epsilon_2}\dots y_{n}^{\epsilon_n}|y_i\in Y\}$ (so $Y^{\emptyset}=\{e\}$ in the empty case). It is easy to see that $\jy=\coprod_{\zeta}Y^{\zeta}$ so that each word $w=y_{1}^{\epsilon_1}y_{2}^{\epsilon_2}\dots y_{n}^{\epsilon_n}$ has a neighborhood base consisting of products $U_{1}^{\epsilon_1}U_{2}^{\epsilon_2}\dots U_{n}^{\epsilon_n}$ where $U_i$ is an open neighborhood of $y_i$ in $Y$. A word $w$ is \textit{reduced} if $y_i=y_{i+1}$ implies $\epsilon_i=\epsilon_{i+1}$ for each $i=1,2,...,n-1$. The empty word is vacuously reduced. The collection of reduced words, of course, forms the free group $F(Y)$ generated by the underlying set of $Y$ and the monoid epimorphism $R:\jy\ra F(Y)$ denotes the usual reduction of words. We can describe the functorial properties of $M_{T}^{\ast}$ using the next definition.
\begin{definition}
A \textbf{topological monoid with continuous involution} is a pair $(M,s)$ where $M$ is a topological monoid with identity $e$ and $s:M\ra M$ is a continuous involution (i.e. $s^2=id$, $s(mn)=s(n)s(m)$, and $s(e)=e$). A morphism $f:(M_1,s_1)\ra(M_2,s_2)$ of two such pairs is a continuous homomorphism $f:M_1\rightarrow M_2$ such that $f$ preserves involution, i.e. $f\circ s_1=s_2\circ f$. Let $\mathscr{TM}^{\ast}$ be the category of topological monoids with continuous involution and continuous, involution-preserving homomorphisms.
\end{definition}

A continuous involution for $\jy$ is given by $(y_{1}^{\epsilon_1}y_{2}^{\epsilon_2}\dots y_{n}^{\epsilon_n})^{-1}=y_{n}^{-\epsilon_n}y_{n-1}^{-\epsilon_{n-1}}\dots y_{1}^{-\epsilon_1}$ making the pair $(\jy,^{-1})$ a topological monoid with continuous involution. One simply uses the pushout property for the coproduct $Y\sqcup Y^{-1}$ and the universal property of free topological monoids to check that
\begin{proposition}
The functor $M_{T}^{\ast}:Top\ra \mathscr{TM}^{\ast}$ is left adjoint to the forgetful functor $\mathscr{TM}^{\ast}\ra Top$.
\end{proposition}
%

We relate free topological monoids to $\pitop$, via the unbased "James map" $u:X\ra \lpsp$, $u(x)(t)=u_x(t)=x\wedge t$. This allows us to define a natural embedding $\mathscr{J}:\jx\ra \lpsp$ taking the empty word to the constant map and $\mathscr{J}(x_{1}^{\epsilon_1}x_{2}^{\epsilon_2}\dots x_{n}^{\epsilon_n})=\ast_{i=1}^{n}\left(u_{x_i}^{\epsilon_i}\right)$. These constructions follow the well known James construction, used, originally by I.M. James, to study the geometry of $\Omega (\Sigma Z, \ast)$ for a connected CW-complex $Z$.
\subsection{The fundamental group $\pi_{1}(\sus)$}

Throughout the rest of this section, we let $P_{X}:X\ra \piztop$ and $P_{\Omega}:\lpsp\ra \pitop$ denote the canonical quotient maps.\\

To compute $\pitop$, we must first understand the algebraic structure of $\pi_{1}(\sus)$. We begin by observing that the James map $u:X\ra \Omega(\Sigma(X_+))$ induces a continuous map $u_{\ast}:\piztop\ra \pi_{0}^{top}(\Omega(\Sigma(X_+)))=\pitop$ on path component spaces. The underlying function $u_{\ast}:\pi_{0}(X)\ra \pi_{1}(\sus)$ induces a group homomorphism $h_{X}:F(\pi_{0}(X))\ra \pi_{1}(\sus)$ on the free group generated by the path components of $X$. In particular, $h_{X}$ takes the reduced word $P_{1}^{\epsilon_1}P_{2}^{\epsilon_2}\dots P_{k}^{\epsilon_k}$ (where $P_i\in \pi_{0}(X)$ and $\epsilon_i\in \{\pm 1\}$) to the homotopy class $[u_{x_1}^{\epsilon_1}\ast u_{x_2}^{\epsilon_2}\ast \dots \ast u_{x_k}^{\epsilon_k}]$ where $x_i\in P_i$ for each $i$. We show that $h_X$ is a group isomorphism. To do this, we make use of the next definition which will also be used later on in section 2.
\begin{definition}
A loop $\alpha\in M(I,\{0,1\};Y,\{y\})$ is \textbf{simple} if $\alpha^{-1}(y)=\{0,1\}$. The subspace of $M(I,\{0,1\};Y,\{y\})$ consisting of simple loops is denoted $\Omega_{s}(Y)$.
\end{definition}
\begin{remark}
$\Omega_{s}$ is not a functor since it is not well defined on morphisms. It is easy to see, however, that $\Omega_{s}(\Sigma((-)_{+})):Top\ra Top$ is a functor.
\end{remark}

The map $X\ra \{\ast\}$ collapsing $X$ to a point induces a retraction $r:\sus\ra \Sigma S^0\cong S^1$. This, in turn, induces a retraction $r_{\ast}:\pitop\ra \pi_{1}^{top}(S^1)\cong \mathbb{Z}$ onto the discrete group of integers. If $\alpha \in \Omega_{s}(\sus)$, then $r\circ \alpha:I\ra S^1$ is a simple loop in $S^1$. But the homotopy class of a simple loop in $S^1$ is either the identity or a generator of $\pi_{1}^{top}(S^1)$. Therefore $r_{\ast}([\alpha])$ must take on the value $1,0$ or $-1$.
\begin{definition}
A simple loop $\alpha\in \Omega_{s}(\sus)$ has \textbf{positive (resp. negative) orienation} if $[\alpha]\in r_{\ast}^{-1}(1)$ (resp. $[\alpha]\in r_{\ast}^{-1}(-1)$). If $[\alpha]\in r_{\ast}^{-1}(0)$, then we say $\alpha$ has no orientation and is \textbf{trivial}. The subspaces of $\Omega_{s}(\sus)$ consisting of simple loops with positive, negative, and no orienation are denoted $\slpsp$, $\Omega_{-s}(\sus)$, and $\Omega_{0s}(\sus)$ respectively.
\end{definition}

The fact that $\mathbb{Z}$ is discrete, allows us to write the loop space $\lpsp$ as the disjoint union $\coprod_{n\in \mathbb{Z}}P_{\Omega}^{-1}(r_{\ast}^{-1}(n))$. Consquently, we may write $\Omega_{s}(\sus)$ as the disjoint union $$\Omega_{s}(\sus)=\slpsp\sqcup\Omega_{0s}(\sus)\sqcup \Omega_{-s}(\sus).$$We also note that $\Omega_{-s}(\sus)=\slpsp^{-1}$. Thus loop inversion give a homeomorphism $\slpsp\cong \Omega_{-s}(\sus)$. The next two lemmas are required to prove the surjectivity of $h_X$.
\begin{lemma} \label{trivialloops}
A simple loop $\alpha\in \Omega_{s}(\sus)$ is null-homotopic if and only if it is trivial.
\end{lemma}
\begin{proof} By definition, a simple loop which has orientation is not null-homotopic. Therefore, it suffices to show that any trivial loop is null-homotopic. If $\alpha$ is trivial, then $\alpha$ does not traverse any hoop of $\sus$, i.e. there is a $t\in (0,1)$ such that $\alpha$ has image in either $X\wedge [0,t]$ or $X\wedge [t,1]$. By 4. of Remark \ref{susfacts}, $\alpha$ is null-homotopic.\end{proof}

We note that the subspaces $P\wedge (0,1)$ for path components $P\in \piz$ are precisely the path components of $X\wedge (0,1)$. Therefore, if $p:I\ra \sus$ is a path such that $p(0)\in P_1\wedge (0,1)$ and $p(1)\in P_2\wedge (0,1)$ for distinct $P_1,P_2\in \pi_{0}(X)$ (i.e. the endpoints of $p$ lie in distinct hoops and are not the basepoint $x_0$), then there is a $t\in (0,1)$ such that $p(t)=x_0$. This implies that the image of each simple loop lies entirely within a single hoop. 
\begin{lemma} \label{orientedsimpleloops}
If simple loops $\alpha$ and $\beta$ have the same orientation and have image in the same hoop $P\wedge I\subseteq \sus$, then they are homotopic.
\end{lemma}
\begin{proof} Suppose $\alpha$ and $\beta$ have positive orientation and image in $P\wedge I$. Since $P\wedge (1,0)$ is a path component of $X \wedge (0,1)$, we may find a $t\in (0,1)$ and a path $p:I\ra X \wedge (0,1)$ such that $p(0)=\alpha(t)$ and $p(1)=\beta(t)$. Now $$\alpha_{[0,t]}\ast p\ast \beta_{[0,t]}^{-1} \text{ and } \beta_{[t,1]}^{-1}\ast p^{-1}\ast \alpha_{[t,1]}$$ are trivial simple loops which by the previous lemma must be null homotopic. This gives fixed endpoint homotopies of paths$$\alpha_{[0,t]}\simeq \beta_{[0,t]}\ast p^{-1}\text{ and }\alpha_{[t,1]}\simeq p\ast \beta_{[t,1]}$$The concatenation of these two gives $$\alpha\simeq \alpha_{[0,t]}\ast\alpha_{[t,1]}\simeq \beta_{[0,t]}\ast p^{-1} \ast p\ast \beta_{[t,1]}\simeq \beta_{[0,t]}\ast \beta_{[t,1]}\simeq \beta.$$One may simply invert loops to prove the case of negative orientation.\end{proof}

We require the next lemma and remark to prove the injectivity of $h_X$.
\begin{lemma} \label{sumcase}
If $w=P_{1}^{\epsilon_1}\dots P_{n}^{\epsilon_n}\in F(\piz)$ is a non-empty reduced word such that $\sum_{i=1}^{n}\epsilon_i\neq 0$, then $h_X(w)$ is not the identity of $\pi_{1}(\sus)$.
\end{lemma}
\begin{proof} The retraction $r:\sus\ra S^1$ induces an epimorphism $r_{\ast}:\pi_{1}(\sus)\ra \mathbb{Z}$ on fundamental groups, where $r_{\ast}([u_x]^{\epsilon})=\epsilon$ for each $x\in X$ and $\epsilon\in \{\pm 1\}$. Therefore, if $\sum_{i=1}^{n}\epsilon_i\neq 0$, then $r_{\ast}(h_{X}(w))=r_{\ast}([u_{x_1}^{\epsilon_1}\ast \dots \ast u_{x_n}^{\epsilon_n}])=\sum_{i=1}^{n}\epsilon_i\neq 0$ (where $x_i\in P_i$) and $h_{X}(w)$ cannot be the identity of $\pi_{1}(\sus)$.\end{proof}
\begin{remark} \label{injectivityremark}\text{ }\\
(1) If $P_{1}^{\epsilon_1}\dots P_{n}^{\epsilon_n}\in F(\piz)$ is a reduced word and $1\leq k\leq m\leq n$, the subword $P_{k}^{\epsilon_k}\dots P_{m}^{\epsilon_m}$ is also reduced.\\
(2) If $P_{1}^{\epsilon_1}\dots P_{n}^{\epsilon_n}\in F(\piz)$ is a reduced word such that $n\geq 2$ and $\sum_{i=1}^{n}\epsilon_i=0$, then there are $i_0,i_1\in\{1,2,...,n\}$ such that $P_{i_0}\neq P_{i_1}$.
\end{remark}
\begin{theorem} \label{fundgrpcomputation}
$h_{X}:F(\piz)\ra \pi_{1}(\sus)$ is an isomorphism of groups.
\end{theorem}
\begin{proof} To show that $h_X$ is surjective, we suppose $\alpha\in \lpsp$ is an arbitrary loop. The pullback $\alpha^{-1}(\sus-\{x_0\})=\coprod_{m\in M}(c_m,d_m)$ is an open subset of $(0,1)$. Each restriction $\alpha_m=\alpha_{[c_m,d_m]}$ is a simple loop, and by \ref{susfacts} (5), all but finitely many of the $\alpha_m$ have image in the simply connected neighborhood $X\wedge [0,\frac{1}{8})\sqcup (\frac{7}{8},1]$. Therefore $\alpha$ is homotopic to a finite concatenation of simple loops $\alpha_{m_1}\ast\alpha_{m_2}\ast \dots \ast \alpha_{m_{n}}$. By Lemma \ref{trivialloops}, we may suppose that each $\alpha_{m_i}$ has orientation $\epsilon_i\in \{\pm 1\}$ and image in hoop $P_i\wedge I$. Lemma \ref{orientedsimpleloops} then gives that $\alpha_{m_i}\simeq u_{x_i}^{\epsilon_i}$ for any $x_i\in P_i$. But then $$h_X(P_{1}^{\epsilon_1}P_{2}^{\epsilon_2}\dots P_{n}^{\epsilon_n})=[u_{x_1}^{\epsilon_1}\ast u_{x_2}^{\epsilon_2}\ast \dots\ast u_{x_n}^{\epsilon_n}]=[\alpha_{m_1}\ast\alpha_{m_2}\ast \dots \ast \alpha_{m_{n}}]=[\alpha].$$
For injectivity, we suppose $w=P_{1}^{\epsilon_1}P_{2}^{\epsilon_2}\dots P_{n}^{\epsilon_n}$ is a non-empty reduced word in $F(\pi_{0}(X))$. It suffices to show that $h_X(w)=[u_{x_1}^{\epsilon_1}\ast u_{x_2}^{\epsilon_2}\ast \dots\ast u_{x_n}^{\epsilon_n}]$ is non-trivial when $x_i\in P_i$ for each $i$. We proceed by induction on $n$ and note that Lemma \ref{sumcase} gives the first step of induction $n=1$. Suppose $n\geq 2$ and $h_X(v)$ is non-trivial for all reduced words $v=Q_{1}^{\delta_1}Q_{2}^{\delta_2}\dots Q_{j}^{\delta_j}$ of length $j<n$. By Lemma \ref{sumcase} it suffices to show that $h_{X}(w)=[u_{x_1}^{\epsilon_1}\ast u_{x_1}^{\epsilon_1}\ast \dots\ast u_{x_n}^{\epsilon_n}]$ is non-trivial when $\sum_{i=1}^{n}\epsilon_i=0$. We suppose otherwise, i.e. that there is a homotopy of based loops $H:I^2\ra \sus$ such that $H(t,0)=x_0$ and $H(t,1)=\left(u_{x_1}^{\epsilon_1}\ast u_{x_2}^{\epsilon_2}\ast\dots\ast u_{x_n}^{\epsilon_n}\right)(t)$ for all $t\in I$. For $j=0,1,...,2n$, we let $b_j=\frac{j}{2n}\times 1\in I^2$. Remark \ref{susfacts} (5) indicates that the singleton $\{x_0\}$ is closed in $\sus$ so that $H^{-1}(x_0)$ is a compact subset of $I^2$. Since each $u_{x_i}^{\epsilon_i}$ is simple we have that $$H^{-1}(x_0)\cap \partial (I^2)=  \{0,1\}\times I\cup I\times \{0\}\sqcup \coprod_{i=1}^{n-1}\{b_{2i}\}$$where $\partial$ denotes boundary in $\mathbb{R}^2$. We also have that $H(b_{2i-1})=u_{x_{i}}^{\epsilon_i}(\frac{1}{2})= x_{i} \wedge \frac{1}{2}\neq x_0$ for each $i=1,...,n$. This allows us to find an $r_0>0$ so that when $U_i=B(b_{2i-1},r_0)\cap I^2$ is the ball of radius $r_0$ about $b_{2i-1}$ in $I^2$, we have $H^{-1}(x_0)\cap \bigcup_{i=1}^{n}U_j=\emptyset$. Now we find an $r_1\in (0,r_0)$ and cover $H^{-1}(x_0)$ with finitely many open balls $V_l=B(z_l,r_1)\cap I^2$ so that $$\left(\bigcup_{l} \overline{V_l}\right)\cap\left(\bigcup_{i=1}^{n}\overline{U_i}\right)=\emptyset\text{ and }H\left(\bigcup_{l} V_l\right)\subseteq [0,\frac{1}{8})\sqcup (\frac{7}{8},1]\wedge X$$(which is possible since we are assuming $H$ to be continuous). Note that if $q:I\ra \bigcup_{l} V_l$ is a path with endpoints $q(0),q(1)\in H^{-1}(x_0)$, then the loop $H\circ q:I\ra \sus$ is based at $x_0$, has image in the simply connected neighborhood $X \wedge [0,\frac{1}{8})\sqcup (\frac{7}{8},1]$, and therefore must be null-homotopic. We pause here to prove the following claim:
\begin{claim} There is no path $q:I\ra \bigcup_{l} V_l$ such that $q(0)=b_{2k}$, $q(1)=b_{2m}$ for distinct $k,m\in\{1,...,n\}$ \end{claim}
\begin{proof} Suppose $q:I\ra \bigcup_{l} V_l$ is such a path and $1\leq k<m\leq n$. But this precisely means that the concatenation $ u_{x_{k+1}}^{\epsilon_{k+1}}\ast u_{x_{k+2}}^{\epsilon_{k+2}}\ast\dots \ast u_{x_{m}}^{\epsilon_m}$ is null-homotopic since $H\circ q$ is null-homotopic and $\left(H\circ q\right)\simeq u_{x_{k+1}}^{\epsilon_{k+1}}\ast u_{x_{k+2}}^{\epsilon_{k+2}}\ast\dots \ast u_{x_m}^{\epsilon_m}$. This means that $h_{X}(P_{k+1}^{\epsilon_{k+1}}P_{k+2}^{\epsilon_{k+2}}\dots  P_{m}^{\epsilon_m})=[u_{x_{k+1}}^{\epsilon_{k+1}}\ast u_{x_{k+2}}^{\epsilon_{k+2}}\ast\dots \ast u_{x_m}^{\epsilon_m}]$ is the identity of $\pi_{1}(\sus)$. But by Remark \ref{injectivityremark} (1) $P_{k+1}^{\epsilon_{k+1}}P_{k+2}^{\epsilon_{k+2}}\dots  P_{m}^{\epsilon_m}$ is a reduced word of length $<n$ and so by our induction hypothesis $h_X(P_{k+1}^{\epsilon_{k+1}}P_{k+2}^{\epsilon_{k+2}}\dots  P_{m}^{\epsilon_m})$ cannot be the identity. This contradiction proves the claim.\end{proof}
\noindent\textit{Proof of Theorem \ref{fundgrpcomputation} cont}. This claim implies that $b_{2i}$ lies in a distinct path component (and consequently connected component) of $\bigcup_{l} V_l$ for each $i=1,...,n$. We let $C_i=\bigcup_{m=1}^{M_i}V_{l_m}^{i}$ be the (path) component of $\bigcup_{l} V_l$ containing $b_{2i}$. But this means that the $b_{2i-1}$, $i=1,...,n$ all lie in the same path component of $I^2-\bigcup_{l} V_l$. Specifically, the subspace $$\left(\partial(I^2)-\bigcup_{l} V_l\right)\cup \left( \partial\left(\bigcup_{i=1}^{n} C_i\right)-\partial(I^2)\right)$$ is path connected and contains each of the $b_{2i-1}$. Since we were able to assume that $\sum_{i}\epsilon_i=0$, we know by Remark \ref{injectivityremark} (2) that there are $i_0,i_1\in \{1,...,n\}$ such that $P_{i_0}\neq P_{i_1}$. We have shown that there is a path $p:I\ra I^2-\bigcup_{l} V_l$ with $p(0)=b_{2i_0-1}$ and $p(1)=b_{2i_1-1}$. But then $H\circ p:I\ra \sus$ is a path with $x_0\notin H\circ p(I)$, $H(p(0))=\frac{1}{2}\wedge x_{i_0}$, and $H(p(1))=\frac{1}{2}\wedge x_{i_1}$. But this is impossible as $H(p(0))$ and $H(p(1))$ lie in different hoops of $\sus$. Therefore $u_{x_1}^{\epsilon_1}\ast u_{x_2}^{\epsilon_2}\ast\dots \ast u_{x_n}^{\epsilon_n}$ cannot be null-homotopic.\end{proof}
\begin{corollary} \label{equalfibers}
The fibers of the map $P_{\Omega}\circ \mathscr{J}:\jx\ra\pitop$ are equal to those of $R\circ \mtpx:\jx\ra \jtop\ra F(\piz)$.
\end{corollary}

Since $u_x\simeq u_y$ if and only if $x$ and $y$ lie in the same path component of $X$, we denote the homotopy class of $u_x$ by $[u_P]$ where $P$ is the path component of $x$ in $X$. Thus $\{[u_P]|P\in\piz \}$ freely generates $\pi_{1}(\sus)$. This computation also indicates that the canonical map $u_{\ast}:\piztop\ra \pitop$ is an injection.
\subsection{Three topologies on $F(\piz)$}

In the effort to recognize the topological structure of $\pitop$, we observe three comparable but potentially distinct topologies on $F(\piz)$. We proceed from the coarsest to the finest topology.
\subsubsection{Free topological groups}
Free topological groups have been extensively studied since their introduction in 1941 by Markov \cite{Markov}. In particular, the free (Markov) topological group $\fty$ on a space $Y$ is the unique (up to topological isomorphism) topological group such that
\begin{enumerate}
\item The canonical injection $\sigma:Y\ra \fty$ is continuous.
\item For each map $f:Y\ra G$ to a topological group $G$, there is a unique continuous homomorphism $\tilde{f}:\fty\ra G$ such that $\tilde{f}\circ \sigma =f$.
\end{enumerate}
$\fty$ is known to exist for every space $Y$ \cite{Kakutani,Samuel} and the algebraic group underlying $\fty$ is the free group $F(Y)$ on the underlying set of $Y$. Most current authors restrict themselves to the case when $Y$ is Tychonoff due to the fact that $Y$ is Tychonoff if and only if $\fty$ is Hausdorff and $\sigma:Y\ra \fty$ is an embedding \cite{Thomas}. In the previous section it was observed that that $u_{\ast}:\piztop\ra \pitop$ is continuous, injective, and has image freely generating the algebraic group $\pi_{1}(\sus)$. Moreover, the composite $h_{X}^{-1}\circ u_{\ast}:\piztop\ra F(\piz)$ is the canonical injection $\sigma:\piztop\ra F(\piz)$. But the topology of $\ftpi$ is the finest group topology on $F(\piz)$ such that $\sigma$ is continuous.
\begin{remark} \label{tgifandonlyif1}
If $\pitop$ is a topological group, then the group isomorphism $h_{x}^{-1}:\ftpi\ra \pitop$ is continuous.
\end{remark}
In section 3, we apply known properties of free topological groups to characterize some topological properties of $\pitop$. The next lemma is well accepted; the abelian version is proved in \cite{Thomas}.
\begin{lemma} \label{presrvingquotients}
$F_{M}$ preserves quotient maps.
\end{lemma}
\begin{proof} Suppose $q:X\ra Y$ is a quotient map and $\sigma_{X}:X\ra \ftx$ and $\sigma_Y:Y\ra \fty$ are the canonical injections. Let $G$ be a topological group and $\tilde{f}:\ftx\ra G$ be a continuous homomorphism such that $\tilde{f}(\ker F_{M}(q))=0$. Suppose $x,x'\in X$ such that $q(x)=q(x')$. Since $F_{M}(q)(\sigma_{X}(x))=\sigma_{Y}(q(x))=\sigma_{Y}(q(x'))=F_{M}(q)(\sigma_{X}(x'))$ and $\tilde{f}$ is constant on the fibers of $F_{M}(q)$, it follows that $f=\tilde{f}\circ \sigma_{X}:X\ra G$ is constant on the the fibers of $q$. This induces a map $k:Y\ra G$ such that $k\circ q=f$. But then $k$ induces a continuous homomorphism $\tilde{k}:\fty\ra G$ such that $\tilde{k}\circ \sigma_{Y}=k$. But then $\tilde{f}\circ \sigma_{X}=f=\tilde{k}\circ \sigma_{Y}\circ q=\tilde{k}\circ F_{M}(q)\circ \sigma_{X}$ and the uniqueness of $\tilde{f}$ gives that $\tilde{f}=\tilde{k}\circ F_{M}(q)$.\end{proof} 
\subsubsection{The reduction topology}
Since we expect $P_{\Omega}:\lpsp\ra\pitop$ to make identifications similar to those made by reduction of words in $\jtop$, a natural choice of topology on the free group $F(\piz)$ is the quotient topology with respect to the reduction homomorphism $R:\jtop\ra F(\piz)$. We may do this in more generality, giving the free group $F(Y)$ the quotient topology with respect to $R:\jy\ra F(Y)$ for each space $Y$. We refer to this topology as the \textit{reduction topology} and denote the resulting group with topology (which will be studied in much greater detail in section 3) as $F_{R}(Y)$. Using arguments similar to those in the proof of Lemma \ref{pi1topfunctorality} the functorality of $F_{R}$ is simple.
\begin{lemma} \label{frfunctorality}
$F_{R}:Top\ra \mathscr{QTG}$ is a functor.
\end{lemma}
\begin{proposition}
The identity homomorphism $\phi_{Y}:\fry\ra\fty$ is continuous. Moreover, $\phi_Y$ is a topological isomorphism if and only if $\fry$ is a topological group.
\end{proposition}
\begin{proof} Let $\sigma ':Y\hookrightarrow \jy$ and $R\circ \sigma ':Y\ra \fry$ be the canonical continuous injections. Since $\fty$ is a topological monoid with continuous involution and $\sigma:Y\ra \fty$ is continuous, there is a unique, involution-preserving, continuous monoid homomorphism $\tilde{\sigma}:\jy\ra \fty$ such that $\sigma=\tilde{\sigma}\circ \sigma '$. But the fibers of $\tilde{\sigma}$ are equal to those of the quotient map $R:\jy \ra \fry$. Therefore the identity $\phi_{Y}:\fry\ra\fty$ is continuous. One direction of the second statement is obvious. If $\fry$ is a topological group, the canonical map $R\circ \sigma ':Y \ra \fry$ induces the continuous identity homomorphism $\phi_{Y}^{-1}:\fty\ra \fry$.\end{proof}

For simplicity of notation, when $Y=\piztop$, we denote the identity as $\Phi_{X}:\frpi\ra \ftpi$. At first glance, the reduction topology appears to be an appropriate choice and its simple description certainly makes it appealing. It is true under certain circumstances that $\pitop\cong \frpi$, in fact, we will fully characterize when this occurs for Hausdorff $X$. To describe the topological isomorphism class of $\pitop$ for arbitrary $X$, however, we make use of a finer topology on $F(\pi_{0}(X))$.
\subsubsection{Semitopological monoids and a quotient reduction topology}
Though a general construction is possible, we restrict ourselves to placing our last topology on the free group $F(\pi_{0}(X))$ as opposed to $\fty$ and $\fry$ which are given for arbitrary $Y$. We make use the two equivalent defects of the category of spaces:
\begin{enumerate}
\item $\pi_{0}^{top}$ fails to preserve finite products.
\item The powers of $P_X:X\ra \piztop$ need not be quotient maps.
\end{enumerate}
\begin{remark}
For a space $X$ and integer $n>1$, the projection maps $pr_{i}:X^n\ra X$ induce maps $(pr_i)_{\ast}:\pi_{0}^{top}(X^n)\ra \pi_{0}^{top}(X)$. Together, these induce a continuous bijection $\psi_n:\pi_{0}^{top}(X^n)\ra \pi_{0}^{top}(X)^n$ which satisfies $\psi_{n}\circ P_{X_n}=P_{X}^{n}$. This confirms the equivalence of two topological defects above. Since $$\pi_{0}^{top}(\jx)\cong \coprod_{n\geq 0}\pi_{0}^{top}((X\sqcup X^{-1})^{n})\text{  and  }\jtop= \coprod_{n\geq 0}(\pi_{0}^{top}(X)\sqcup\pi_{0}^{top}(X)^{-1})^{n}$$we have, for each space $X$, a canonical continuous bijection $\psi_{X}:\pi_{0}^{top}(\jx)\ra \jtop$. The quotient maps $P_X:X\ra \piztop$ and $P_{\jx}:\jx\ra \pi_{0}^{top}(\jx)$ satisfy $\psi_{X}\circ P_{\jx}=\mtpx$.
\end{remark}

Now $\pi_{0}^{top}(\jx)$ inherits the monoid structure of $\jtop$, however, $\psi_{X}$ need not be open. In fact, we can fully characterize when this occurs.
\begin{corollary} \label{monoidequivalentst}
The following are equivalent:
\begin{enumerate}
\item $\pi_{0}^{top}(\jx)$ is a topological monoid.
\item $\psi_{X}$ is a topological isomorphism.
\item $\mtpx$ is a quotient map.
\item For each integer $n\geq 2$, $(P_X)^{n}:X^n\ra \piztop^{n}$ is a quotient map.
\end{enumerate}
\end{corollary}
\begin{proof} 2. $\Rightarrow$ 1. is obvious. 1. $\Rightarrow$ 2. If $\pi_{0}^{top}(\jx)$ is a topological monoid, the canonical inclusion $\piztop\hookrightarrow \pi_{0}^{top}(\jx)$ induces a continuous homomorphism $\jtop\ra \pi_{0}^{top}(\jx)$ which is the inverse of $\psi_{X}$. 2. $\Leftrightarrow$ 3. follows from the fact that $\psi_{X}\circ P_{\jx}=\mtpx$ and that $P_{\jx}$ is a quotient map. 3. $\Leftrightarrow$ 4. follows from the definition $\mtpx=\coprod_{n\geq 0}(P_X\sqcup P_X)^{n}$.\end{proof}

We give the next definition since multiplication in $\pi_{0}^{top}(\jx)$ is not always continuous.
\begin{definition} \emph{
A \textbf{semitopological monoid} is a monoid $M$ with topology such that multiplication $M\times M\ra M$ is separately continuous. A \textbf{semitopological monoid with continuous involution} is a pair $(M,s)$ where $M$ is a semitopological monoid and $s:M\ra M$ a continuous involution. A morphism of semitopological monoids with continuous involution is a continuous, involution-preserving, monoid homomorphism. The category described here will be denoted $\mathscr{STM}^{\ast}$.}
\end{definition}

The proof of the functorality of $\pi_{0}^{top}(\jx)$ is similar to that of Lemmas \ref{pi1topfunctorality} and \ref{frfunctorality}.
\begin{lemma}
$\pi_{0}^{top}(M_{T}^{\ast}(-)):Top\ra \mathscr{STM}^{\ast}$ is a functor. Since every topological monoid is a semitopological monoid, the continuous monoid isomorphisms $\psi_{X}:\pi_{0}^{top}(\jx)\ra \jtop$ are components of a natural transformation.
\end{lemma}

We are now interested in the composite monoid epimorphism:$$\xymatrix{\pi_{0}^{top}(\jx) \ar[r]^{\psi_{X}}  &  \jtop  \ar[r]^{R} & F(\pi_{0}(X))}.$$We give $F(\pi_{0}(X))$ the quotient topology with respect to $R\circ \psi_{X}$ and denote the resulting group with topology as $\fotop$. Functorality follows as before and the universal property of quotient spaces gives an immediate relationship to $\frpi$.
\begin{lemma} \label{finertopology}
$X\mapsto \fotop$ gives a functor $Top\ra \mathscr{QTG}$.
\end{lemma}

To relate $\fotop$ and $\frpi$ we observe that the vertical arrows in the commuting diagram$$\xymatrix{ \pi_{0}^{top}(\jx) \ar[r]^{\psi_{X}} \ar[d]_{R\circ \psi_{X}} & \jtop \ar[d]^{R}\\ \fotop \ar[r]_{\Psi_{X}} & \frpi }$$are quotient. Consequently the identity homomorphism $\Psi_{X}:\fotop\ra \frpi$ is always continuous.
\subsubsection{A comparison of $\fotop$, $\frpi$, and $\ftpi$}
We begin our comparison of the three topologies on $F(\pi_{0}(X))$ using the commutative diagram$$\xymatrix{ & \piztop \ar[dl]_{\sigma} \ar[d]^{\sigma} \ar[dr]^{\sigma} \\
\fotop \ar[r]_{\Psi_{X}} & \frpi \ar[r]_{\Phi_{X}} & \ftpi }$$Here $\sigma:\piztop\ra F(\pi_{0}(X))$ is the canonical injection which is continuous with respect to all three topologies. The horizontal maps are the continuous identity homomorphisms.\\

The next diagram gives a more detailed comparison of the topological structures.
$$\xymatrix{ & \jx \ar[dl]_{P_{\jx}} \ar[r]^{R} \ar[d]^{\mtpx} & \frx \ar[ddl]^{F_{R}(P_X)} \ar[d]^{\phi_{X}} \\
\pi_{0}^{top}(\jx) \ar[r]_{\psi_{X}} \ar[d]_{R\circ \psi_{X}} & \jtop \ar[d]_{R} & \ftx \ar[d]^{F_{M}(P_X)} \\
\fotop \ar[r]_{\Psi_{X}} & \frpi \ar[r]_{\Phi_{X}} & \ftpi }$$The maps $\psi_{X},\phi_{X},\Psi_{X},\Phi_{X}$ are continuous identity homomorphisms. By Lemma \ref{presrvingquotients}, $F_{M}(P_X)$ is a quotient map. In fact, the only maps other than the identities that need not be quotient are $\mtpx$ and $F_{R}(P_X)$.
\begin{proposition} \label{classify1} $\Psi_{X}:\fotop\ra \frpi$ is a topological isomorphism if and only if $F_{R}(P_X):\frx\ra \frpi$ is a quotient map. If $(P_X)^{n}$ is a quotient map for all $n\geq 2$, then $F_{R}(P_X)$ must be quotient.
\end{proposition}
\begin{proof}
The first statement is clear from the placement of the quotient maps in the second diagram. If $(P_X)^{n}$ is a quotient map for all $n\geq 2$, then $\mtpx:\jx\ra \jtop$ is a quotient map. But since $R\circ \mtpx=F_{R}(P_X)\circ R$ where the $R$ are the appropriate reduction quotient maps $F_{R}(P_X)$ must also be quotient.
\end{proof}
\begin{proposition}  \label{classify3}
The following are equivalent:
\begin{enumerate}
\item $\fotop$ is a topological group.
\item $\Psi_{X}$ and $\Phi_{X}$ are topological isomorphisms.
\item The composite $F_{M}(P_X)\circ \phi_{X}:\frx\ra \ftpi$ is a quotient map.
\end{enumerate}
\end{proposition}
\begin{proof}
2. $\Rightarrow$ 1. is obvious. 1. $\Rightarrow$ 2. If $\fotop$ is a topological group, then the canonical continuous injection $\sigma:\piztop\ra \pi_{0}^{top}(\jx)\ra \fotop$ induces the continuous identity $\ftpi\ra \fotop$. But since the topology of $\fotop$ is always finer than $\frpi$ which is finer than that of $\ftpi$ (from the first diagram), all three must be topologically isomorphic. 2. $\Leftrightarrow$ 3. simply follows from observing the placement of the quotient maps in the second diagram.
\end{proof}
\subsection{Relating $\fotop$ and $\pitop$}

We use the following commutative diagram to relate the quasitopological groups $\pitop$ and $\fotop$.
$$\xymatrix{ \lpsp \ar[d]_{P_{\Omega}} & \jx \ar[l]_{\mathscr{J}} \ar[d]_{P_{\jx}} \ar[dr]^{\mtpx}\\
\pitop & \pi^{top}_{0}(\jx) \ar[l]_{\mathscr{J}_{\ast}} \ar[d]^{R\circ \psi_{X}} \ar[r]_{\psi_{X}}  & \jtop \ar[d]^{R}\\
& \fotop \ar@{-->}[ul]^{h_{X}} \ar[r]_{\Psi_{X}} & \frpi }$$Here $\mathscr{J}_{\ast}=\pi_{0}^{top}(\mathscr{J})$ is the map induced by $\mathscr{J}$ on path component spaces. We recall from Corollary \ref{equalfibers}, that the fibers of the composites $P_{\Omega}\circ \mathscr{J}:\jx\ra \pitop$ and $R\circ \mtpx:\jx\ra \frpi$ are equal. Since $\Psi_{X}$ is the identity homomorphism, it is easy to see that the fibers of the quotient map $(R\circ \psi_{X})\circ P_{\jx}:\jx\ra \fotop$ are equal to the fibers of the map $P_{\Omega}\circ \mathscr{J}:\jx\ra \pitop$. The universal property of quotient spaces, then indicates that the group isomorphism $h_{X}:\fotop\ra \pitop$ of section 2.3 is continuous. Our main result asserts that for arbitrary $X$, $h_{X}$ is also an open map and hence
%
%
%
\begin{theorem}\label{mainresult}
For each space $X$, $h_{X}:\fotop\ra \pitop$ is a topological isomorphism.
\end{theorem}

This result is a "computation" in the sense that it fully identifies the isomorphism class of $\pitop$ in $\mathscr{QTG}$. Specifically, $\pitop$ is the quotient (via reduction) of the path component space of the free topological monoid with continuous involution on $X$. We will prove Theorem \ref{mainresult} in the next section but first we note some immediate consequences. First, Proposition \ref{classify3} allows us to characterize the spaces $X$ for which $\pitop$ is a topological group.
\begin{corollary} \label{tgclassification1}The following are equivalent:
\begin{enumerate}
\item $\pitop$ is a topological group.
\item $\pitop$ is topologically isomorphic to the free topological group $\ftpi$.
\item $\Phi_{X}$ and $\Psi_{X}$ are topological isomorphisms.
\end{enumerate}
\end{corollary}
\begin{proof}
3. $\Rightarrow$ 2. follows directly from \ref{classify3} and 2. $\Rightarrow$ 1. is obvious. 1. $\Rightarrow$ 3. If $\pitop$ is a topological group, the canonical map $u_{\ast}:\piztop\ra \pitop$ induces the continuous homomorphism $h_X:\ftpi\ra \pitop$. Then the identity $id=h_{X}^{-1}\circ h_{X}:\ftpi\ra\pitop\ra  \fotop$ is continuous and the topologies of $\fotop$, $\frpi$, and $\ftpi$ agree. Therefore $\Phi_{X}$ and $\Psi_{X}$ are topological isomorphisms.
\end{proof}
We admit that this characterization is not particularly illuminating. Much of section 3 is devoted to translating the condition that $\Phi_{X}$ and $\Psi_{X}$ be a topological isomorphism into familiar conditions from topology. We do however, have the following sufficient conditions.
\begin{corollary} \label{tgclassification2}
If $\frx$ is a topological group, then so is $\pitop$. 
\end{corollary}
\begin{proof}
If $\frx$ is a topological group, then $\phi_{X}:\frx\ra \ftx$ is a homeomorphism and the composite $F_{M}(P_X)\circ \phi_{X}:\frx\ra \ftpi$ is a quotient map. By Theorem \ref{mainresult} and Proposition \ref{classify3} we have topological isomorphisms $\pitop\cong \fotop\cong \frpi\cong \ftpi$ proving that $\pitop$ is indeed a (free) topological group.
\end{proof}
\begin{example} \label{totallypathdisconnected}
 \emph{
This computation becomes remarkably simple for totally path disconnected spaces (i.e. spaces $X$ such that $P_X:X\cong \piztop$). In this case $(P_X)^n$ is a homeomorphism for all $n\geq 2$ and so $$\pitop\cong \fotop\cong\frpi\cong \frx.$$
}
\end{example}
We may also give a nice characterization of discreteness.
\begin{corollary} \label{discretenessofpi1top}
$\pitop$ is discrete if and only if $\piztop$ is discrete.
\end{corollary}
\begin{proof} If $\piztop$ is discrete, then so is $\jtop$ and its quotient $\frpi$. Since the topology of $\pitop$ is finer than that of $\frpi$, $\pitop$ must also be discrete. Conversely, if $\pitop$ is discrete, the continuity of the injection $u_{\ast}:\piztop\ra \pitop$ gives the discreteness of $\piztop$. \end{proof} 
\subsection{A Proof of Theorem \ref{mainresult}}
To prove that $h_{X}$ is open, we take the approach of factoring the quotient map $P_{\Omega}:\lpsp\ra \pitop$ into the composite of functions, not all of which will be continuous. We begin by studying the topology of simple loops (defined in section 2.2) and decomposing arbitrary loops into words of oriented simple loops.\\
\\
\textit{Step 1: The topology of simple loops}\\

Throughout the rest of this section let $U=X\wedge \left[0,\frac{1}{8}\right)\sqcup\left(\frac{7}{8},1\right]$. This is an arc-connected, contractible neighborhood and by the definition of $\mathscr{B}_{\sus}$ contains all basic open neighborhoods of the basepoint $x_0$. We now prove a basic property of open neighborhoods of simple loops considered as elements in the free path space. Recall that basic open neighborhoods in $M(I,\sus)$ are those described in Lemma \ref{basis} with respect the the basis $\mathscr{B}_{\sus}$.
\begin{lemma} \label{simpleloopneighborhood}
Suppose $0<\epsilon<\frac{1}{8}$ and $W=\bigcap_{i=1}^{m}\langle K_{m}^{i},W_i\rangle$ is a basic open neighborhood of simple loop $\alpha:I\ra \sus$ in the free path space $M(I,\sus)$. There is a basic open neighborhood $V_0$ of $x_0$ in $\sus$ contained in $X\wedge [0,\epsilon)\sqcup (1-\epsilon,1]$ in $\sus$ and a basic open neighborhood $V=\bigcap_{j=1}^{n}\langle K_{n}^{j},V_j\rangle$ of $\alpha$ in $M(I,\sus)$ contained in $W$ such that:
\begin{enumerate}
\item $V_0=V_1=V_2=\dots=V_l=V_k=V_{k+1}=\dots=V_n$ for integers $1\leq l<k\leq n$.
\item The open neighborhoods $V_{l+1},\dots,V_{k-1}$ are of the form $A\wedge (a,b)$ where $A\in\mathscr{B}_{X}$ and $b-a<\epsilon$.
\end{enumerate}
\end{lemma}
\begin{proof} Let $V_0=(W_1\cap W_m)\cap (X\wedge [0,\epsilon)\sqcup (1-\epsilon,1])\subset U$. Since $\mathscr{B}_{\sus}$ is closed under finite intersection $V_0\in \mathscr{B}_{\sus}$. There is an integer $M>3$ such that $m$ divides $M$ and $\alpha(K_{M}^{1}\sqcup K_{M}^{M})\subseteq V_0$. Since $\alpha$ is simple we have $\alpha([\frac{1}{M},\frac{M-1}{M}])\subset X\wedge (0,1)$. When $p=2,...,M-1$ and $K_{M}^{p}\subseteq K_{m}^{i}$ we may cover $\alpha(K_{M}^{p})$ with finitely many open neighborhoods contained in $W_i\cap (X\wedge (0,1))$ of the form $A\wedge (a,b)$ where $A\in \mathscr{B}_{X}$ and $b-a<\epsilon$. We then apply the Lebesgue lemma to take even subdivisions of $I$ to find open neighborhoods $Y_i=\bigcap_{q=1}^{N_p}\langle K_{N_p}^{q},Y_{p,q}^{i}\rangle\subseteq \langle I,W_i\rangle$ of the restricted path $\alpha_{K_{M}^{p}}$. Here each $Y_{p,q}^{i}$ is one of the open neighborhoods $A\wedge (a,b)\subseteq W_i$. We now use the induced neighborhoods of section 1.2 to define$$V=\langle K_{M}^{1}\sqcup K_{M}^{M},V_0\rangle\cap \bigcap_{p=2}^{M-1}\left((Y_{i})^{K_{M}^{p}}\right).$$This is an open neighborhood of $\alpha$ by definition and it suffices to show that $V\subseteq W$. We suppose $\beta\in V$ and show that $\beta(K_{m}^{i})\subseteq W_i$ for each $i$. Clearly, $\beta(K_{M}^{1}\sqcup K_{M}^{M})\subseteq W_1\cap W_m$. If $p=2,...,M-1$ and $K_{M}^{p}\subseteq K_{m}^{i}$, then $\beta_{K_{M}^{p}}\in V_{K_{M}^{p}}\subseteq Y_i$ and $\beta(K_{M}^{p})=\beta_{K_{M}^{p}}(I)\subseteq \bigcup_{q=1}^{N_p}Y_{p,q}^{i}\subseteq W_i$. We may write $V$ as $V=\bigcap_{j=1}^{n}\langle K_{n}^{j},V_j\rangle$ simply by finding an integer $n$ which is divisible by $M$ and every $N_p$ and reindexing the open neighborhoods $V_0$ and $Y_{p,q}^{i}$. In particular, we can set $V_j=V_0$ when $K_{n}^{j}\subseteq K_{M}^{1}\cup K_{M}^{M}$. Additionally, if $H_{K_{M}^{p}}^{-1}:I\ra K_{M}^{p}$ is the unique linear homeomorphism (as in section 1.2) then we let $V_j=Y_{p,q}^{i}$ whenever $$K_{n}^{j}\subseteq H_{K_{M}^{p}}^{-1}\left(K_{N_p}^{q}\right)\subseteq K_{M}^{p}\subseteq K_{m}^{i}.$$It is easy to see that both 1. and 2. in the statement are satisfied by $V$.\end{proof}

We note some additional properties of the neighborhood $V$ constructed in the previous lemma:
\begin{remark}\label{simpnbhdremark1}
\emph{
For each path $\beta\in V$ we have,$$\beta\left(\left(\bigcup_{j=1}^{l}K_{n}^{j}\right)\cup \left(\bigcup_{j=k}^{n}K_{n}^{j}\right)\right)\subseteq U \text{ and }x_0\notin \beta\left(\bigcup_{j=l+1}^{k-1}K_{n}^{j}\right)$$This follows directly from the conditions 1. and 2. in the lemma.
}
\end{remark}
\begin{remark} \label{simpnbhdremark2}
\emph{
It was noted in section 2.2 that there are disjoint open neighborhoods $W_+$, $W_0$, and $W_{-}$ in $M(I,\sus)$ containing $\slpsp$, $\Omega_{0s}(\sus)$, and $\Omega_{-s}(\sus)$ respectively. Consequently, If $\alpha$ has positive orientation then we may take $V\subseteq W_+$ such that $V\cap \Omega_{s}(\sus)\subseteq \slpsp$, i.e. all simple loops in $V$ also have positive orientation. The same holds for the negative and trivial case. In some sense, this means that $V$, when thought of as an instruction set, is "good enough" to distinguish orientations of simple loops.
}
\end{remark}
\begin{remark} \label{nearbysimpleloops} \emph{
We now give a construction, necessary for step 5, which produces a simple loop $\mathscr{S}_{V}(\beta)\in V$ for each path $\beta\in V$. For brevity, we let $[0,r]=\bigcup_{j=1}^{l}K_{n}^{j}$, $[r,s]=\bigcup_{j=l+1}^{k-1}K_{n}^{j}$, and $[s,1]=\bigcup_{j=k}^{n}K_{n}^{j}$. We define $\mathscr{S}_{V}(\beta)$ piecewise by letting $\mathscr{S}_{V}(\beta)$ be equal to $\beta$ on the middle interval $[r,s]$ (i.e. $\mathscr{S}_{V}(\beta)_{[r,s]}=\beta_{[r,s]}$). We then demand that $\mathscr{S}_{V}(\beta)$ restricted to $[0,r]$ is an arc in $V_0$ connecting $x_0$ to $\beta(r)$ and similarly $\mathscr{S}_{V}(\beta)$ restricted to $[s,1]$ is an arc in $V_0$ connecting $\beta(s)$ to $x_0$. Since the image of $\mathscr{S}_{V}(\beta)$ on $[0,r]\cup [s,1]$ remains in $V_0$, it follows that $\mathscr{S}_{V}(\beta)\in V$. Additionally, Remark \ref{simpnbhdremark1} and the use of arcs to define $\mathscr{S}_{V}(\beta)$ implies that $\mathscr{S}_{V}(\beta)$ is a simple loop.}
\end{remark}
\noindent \textit{Step 2: Decomposition of arbitrary loops}\\

Here we assign to each loop in $\sus$, a (possibly empty) word of simple loops with orientation. We again use the observation, that $\slpsp$ and $\Omega_{-s}(\sus)=\slpsp^{-1}$ are disjoint homeomorphic subspaces of $M(I,\sus)$. The free topological monoid on $\slpsp\sqcup \slpsp^{-1}$ is topologically isomorphic to the free topological monoid with continuous involution $\jlpsp$. From here on we make no distinction between the one letter word $\alpha^{-1}$ in $\jlpsp$ and the inverse loop $\beta=\alpha^{-1}\in \slpsp^{-1}$. Similarly, a basic open neighborhood of $\alpha^{-1}$ in $\jlpsp$ corresponds to an open neighborhood of $\beta$ in $\slpsp^{-1}$.\\

We now define a "decomposition" function $\mathscr{D}:\lpsp\ra \jlpsp$. In step 5, we refer to the details of the construction of $\mathscr{D}$.\\
\\
\noindent\textbf{Decomposition:} Suppose $\beta\in M(I,\{0,1\};\sus,\{x_0\})$ is an arbitrary loop. First, if $\beta$ has image contained in $U$ (i.e. $\beta\in \langle I,U\rangle$), then we let $\mathscr{D}(\beta)=e$ be the empty word. Suppose then that $\beta(I)\nsubseteq U$. The pullback $\beta^{-1}(X\wedge (0,1))=\coprod_{m\in M}(c_m,d_m)$ is open in $I$ where $M$ is a countable indexing set with ordering induced by the ordering of $I$. Each restricted loop $\beta_m=\beta_{[c_m,d_m]}:I\ra \sus$ is a simple loop. Remark \ref{susfacts} (4) implies that all but finitely many of these simple loops have image in $U$ and so we may take $m_{1}<...<m_{k}$ to be the indices of $M$ corresponding to those $\beta_{m_{1}},...,\beta_{m_{k}}$ with image not contained in $U$. Note that if $C=I-\bigcup_{i=1}^{k}(c_{m_{i}},d_{m_{i}})$, then we have $\beta\in \langle C,U\rangle$. If none of the $\beta_{m_{i}}$ have orientation, we again let $\mathscr{D}(\beta)=e$. On the other hand, if one of the $\beta_{m_{i}}$ has orientation, we let $m_{i_1}<...<m_{i_n}$ be the indices which correspond to the simple loops $\beta_j=\beta_{m_{i_j}}$ which have either positive or negative orientation. We then let $\mathscr{D}(\beta)$ be the word $\beta_1\beta_2...\beta_n$ in $\jlpsp$.
\begin{remark} \emph{
Informally, $\mathscr{D}(\beta)$ denotes the word composed of the simple loops of $\beta$ which contribute a letter in the unreduced word of the homotopy class $[\beta]$. We may suppose that $\beta_i$ has image in $P_i\wedge I$ and orientation $\epsilon_i\in\{\pm 1\}$, or equivalently that $[\beta_i]=[u_{P_i}]^{\epsilon_i}$. Clearly $\beta\simeq \ast_{i=1}^{n}\beta_i$ and $\varphi_{X}([\beta])=R(P_{1}^{\epsilon_1}P_{n}^{\epsilon_n}\dots P_{n}^{\epsilon_n})\in F(\pi_{0}(X))$.
}
\end{remark}
\noindent \textit{Step 3: An isomorphism lemma}\\

The James map $u:X\ra \lpsp$ has image in $\Omega_{+s}(\sus)$, and the map $u:X\ra \slpsp$ with restricted codomain induces a continuous bijection $u_{\ast}:\piztop\ra \pi_{0}^{top}(\slpsp)$ on path component spaces. The fact that $u_{\ast}$ is also a homeomorphism follow from the argument use to prove the next lemma. For a map $f:X\ra Y$, we let $f_{\ast\ast}=\pi_{0}^{top}(M_{T}^{\ast}(f))$ be the induced continuous involution-preserving monoid homomorphism.
\begin{lemma}
The James map $u:X\ra \slpsp$ induces a natural isomorphism of semitopological monoids with continuous involution $u_{\ast\ast}:\pi_{0}^{top}(\jx)\ra \pi_{0}^{top}\left(\jlpsp\right)$.
\end{lemma}
\begin{proof} We note that on generators $u_{\ast\ast}$ is given by $u_{\ast\ast}(P)=[u_P]$. The naturality of $\psi:\pi_{0}^{top}(M_{T}^{\ast}(-))\ra M_{T}^{\ast}(\pi_{0}^{top}(-))$ applied to the James map makes the following diagram commute in the category of monoids (without topology)$$\xymatrix{ \pi_{0}^{top}(\jx) \ar[d]_{u_{\ast\ast}} \ar[rr]^{\psi_{X}}_{\cong} & & \jtop \ar[d]_{\cong}^{M_{T}^{\ast}(u_{\ast})} \\
\pi_{0}^{top}(\jlpsp) & & M_{T}^{\ast}(\pi_{0}^{top}(\slpsp)) \ar[ll]^{\psi_{\slpsp}^{-1}}_{\cong}
}$$Since $u_{\ast}$ is a bijection, $M_{T}^{\ast}(u_{\ast})$ is a monoid isomorphism. Therefore $u_{\ast\ast}$ is a continuous, involution-preserving monoid isomorphism and it suffices to construct an inverse. Let $r:\slpsp\ra M((0,1),(0,1)\times X)$ be the map taking each positively oriented simple loop $\alpha:I\ra \sus$ to the restricted map $\alpha|_{(0,1)}:(0,1)\ra X\wedge (0,1)\cong X\times (0,1)$ and $p:M((0,1),X\times (0,1))\ra M((0,1),X)$ be post-composition with the projection $X\times (0,1)\ra X$. For any $t\in (0,1)$, we may consider the composite map$$\xymatrix{ v:\slpsp \ar[r]^{j_t} & (0,1)\times \slpsp \ar[rr]^{id\times (p\circ r)} & &  (0,1)\times M((0,1),X) \ar[rr]^{ev}  && X}$$where $j_t(\alpha)=(t,\alpha)$ and $ev$ is the evaluation map $ev(t,f)=f(t)$. If $\alpha\in \slpsp$ such that $\alpha(t)=x\wedge s$, then $v(\alpha)=x$. It is easy to check that the continuous homomorphism $v_{\ast\ast}:\pi_{0}^{top}(\jlpsp) \ra \pi_{0}^{top}(\jx)$ is the inverse of $u_{\ast\ast}$ since on the generator $[u_P]$ of $\pi_{0}^{top}(\jlpsp)$, we have $v_{\ast\ast}([u_P])=v_{\ast\ast}([u_x])=P$ (where $x\in P$).\end{proof}
\noindent\textit{Step 4: Factoring $P_{\Omega}$}\\

We factor the quotient map $P_{\Omega}:\lpsp\ra \pitop$ into a composite using the following functions:
\begin{enumerate}
\item The decomposition function $\mathscr{D}:\lpsp\ra \jlpsp$.
\item The quotient map $P_s:\jlpsp\ra \pi_{0}^{top}(\jlpsp)$ identifying path components (homotopy classes of positively oriented simple loops).
\item The topological isomorphisms $u_{\ast\ast}^{-1}:\pi_{0}^{top}\left(\jlpsp\right)\ra \pi_{0}^{top}(\jx)$ of step 3.
\item The quotient map $R\circ\psi_{X}:\pi_{0}^{top}(\jx)\ra \fotop$.
\item The continuous group isomorphism $h_{X}:\fotop\ra \pitop$
\end{enumerate}
We let $K=R\circ\psi_{X}\circ \rho_{X}\circ P_{s}\circ \mathscr{D}:\lpsp\ra \fotop$ be the composite of 1.-4. and $K'=R\circ(\psi_{X}\circ \rho_{X})\circ P_{s}:\jlpsp\ra \fotop$ be the continuous (and even quotient) composite of 2.-4.
\begin{lemma} \label{diagram}
The following diagram commutes:
$$\xymatrix{
\lpsp \ar[dddr]_{K} \ar[ddd]_{P_{\Omega}} \ar[r]^{\mathscr{D}} & \jlpsp \ar[d]^{P_{s}} \\ 
& \pi_{0}^{top}\left(\jlpsp\right) \ar[d]_{\cong}^{u_{\ast\ast}^{-1}} \\
 & \pi_{0}^{top}(\jx)  \ar[d]^{R\circ\psi_{X}} \\
\pitop & \fotop \ar[l]^{h_X}
}$$
\end{lemma}
We note that the function $\mathscr{D}$ will not be continuous even when $X$ contains only a single point (i.e. $\sus\cong S^1$). This is illustrated by the fact that any open neighborhood $V=\bigcap_{j=1}^{n}\langle K_{n}^{j},U_j\rangle$ of a concatenation $\alpha\ast \alpha^{-1}$ of a simple loop $\alpha$ with it's inverse contains a simple $\beta$ which has no orientation. $\beta$ may be found by "pulling" the middle of $\alpha\ast\alpha^{-1}$ off of $x_0$ within a sufficiently small neighborhood of $x_0$. For this reason, our main difficulty lies in proving the continuity of $K$.\\
\\
\textit{Step 5: Continuity of $K$}
\begin{lemma} \label{continuouscomposite}
The composite $K:\lpsp\ra \fotop$ is continuous.
\end{lemma}
\begin{proof} Suppose $W$ is open in $\fotop$ and $\beta\in K^{-1}(W)$. We now refer to the details of the decomoposition of $\beta$ in step 2. If $\beta$ has image in $U$, then clearly $\beta\in \langle I,U\rangle\subseteq \mathscr{D}^{-1}(e)\subseteq K^{-1}(W)$. Suppose, on the other hand, that some simple loop restriction $\beta_{m_i}$ has image intersecting $\sus-U$ and $\mathscr{D}(\beta)=\beta_{1}\beta_{2}...\beta_{n}$ is the (possibly empty) decomposition of $\beta$. Recall from our the notation in step 2, that $\beta_j=\beta_{m_{i_j}}$, $j=1,...,n$ are the $\beta_{m_i}$ with orientation. Since $K'=R\circ\psi_{X}\circ \rho_{X}\circ P_{s}$ is continuous, $(K')^{-1}(W)$ is an open neighborhood of $\mathscr{D}(\beta)$ in $\jlpsp$.\\
\\
We recall that $\beta\in \langle C,U\rangle$, where $C=I-\bigcup_{i=1}^{k}(c_{m_i},d_{m_i})$. We construct the rest of desired open neighborood of $\beta$ by defining an open neighborhood of each $\beta_{m_i}$ and taking the intersection of the induced neighborhoods.\\ 
\\
If $i\neq i_j$ for any $j=1,...,n$, then $\beta_{m_i}$ does not appear as a letter in the decomposition of $\beta$ and must be trivial. We apply Lemma \ref{simpleloopneighborhood}, to find an open neighborhood $V_{i}=\bigcap_{l=1}^{N_i}\langle K_{N_i}^{l},V_{l}^{i}\rangle$ of $\beta_{m_i}$ in $M(I,\sus)$ which satisfies both 1. and 2. in the statement. By Remark \ref{simpnbhdremark2}, we may also assume that $V_i\cap \Omega_{s}(\sus)\subseteq \Omega_{0s}(\sus)$, i.e all simple loops in $V_i$ are trivial.\\
\\
If $i=i_j$ for some $j=1,...,n$, then $\beta_{j}=\beta_{m_{i_j}}$ has orienation $\epsilon_j$. Since $(K')^{-1}(W)$ is an open neighborhood of $\mathscr{D}(\beta)=\beta_{1}\beta_{2}...\beta_{n}$ in $\jlpsp$ and basic open neighborhoods in $\jlpsp$ are products of open neighborhoods in $\slpsp$ and $\slpsp^{-1}$, we can find basic open neighborhoods $V_{i_j}=\bigcap_{l=1}^{N_{i_j}}\langle K_{N_{i_j}}^{l},U_{l}^{i_j}\rangle$ of $\beta_{j}$ in $M(I,\sus)$ such that $$W_{j}=V_{i_j}\cap\slpsp^{\epsilon_j}\text{ and }\beta_{1}\beta_{2}...\beta_{n}\in W_{1}W_{2}...W_{n}\subseteq (K')^{-1}(W).$$We assume each $V_{i_j}$ satisfies 1. and 2. of Lemma \ref{simpleloopneighborhood} and by Remark \ref{simpnbhdremark2} that $V_{i_j}\cap \Omega_{s}(\sus)\subseteq \lpsp^{\epsilon_{j}}$.\\
\\
Let $$\mathscr{U}=\langle C,U\rangle \cap\left(\bigcap_{i=1}^{k}\left(V_{i}^{[c_{m_i},d_{m_i}]}\right)\right)$$so that $\mathscr{V}=\mathscr{U}\cap \lpsp$ is an open neighborhood of $\beta$ in the loop space. We claim that each loop $\gamma\in \mathscr{V}$ is homotopic to a loop $\gamma '$ such that $\mathscr{D}(\gamma ')\in  (K')^{-1}(W)$. If this is done, we have $$h_{X}(K(\gamma))=P_{\Omega}(\gamma)=P_{\Omega}(\gamma ')=h_{X}(K(\gamma '))$$and since $h_{X}$ is a bijection, $$K(\gamma)=K(\gamma ')=K'(\mathscr{D}(\gamma '))\in W$$This gives $K(\mathscr{V})\subseteq W$, proving the continuity of $K$.\\
\\
We define $\gamma '$ piecewise and begin by setting $\gamma '(C)=x_0$. The restricted path $\gamma_i=\gamma_{[c_{m_i},d_{m_i}]}:I\ra \sus$ lies in the open neighborhood $\mathscr{U}_{[c_{m_i},d_{m_i}]}\subseteq V_i$. We now define $\gamma '$ on $[c_{m_i},d_{m_i}]$ by using the construction of Remark \ref{nearbysimpleloops}. We set $${\gamma_i}'=(\gamma ')_{[c_{m_i},d_{m_i}]}=\mathscr{S}_{V_i}(\gamma_{i})$$which by construction is a simple loops in $V_i$. Intuitively, we have replaced the portions of $\gamma$ which are close to $x_0$ ("close" meaning with respect to $\mathscr{V}$) with arcs and constant paths. Since ${\gamma_i}'=\gamma ' (C)=x_0\in U$ and $(\gamma ')_{[c_{m_i},d_{m_i}]}\in V_i$ for each $i$, it follows that $\gamma '\in \mathscr{V}$. Moreover, since $\gamma(t)\neq \gamma '(t)$ only when $\gamma(t)$ and $\gamma '(t)$ both lie in the path connected, contractable neighborhood $U$, it is obvious that $\gamma\simeq \gamma '$. It now suffices to show that $\mathscr{D}(\gamma ')\in  (K')^{-1}(W)$. We begin by checking which of the simple loops ${\gamma_i}'\in V_i\cap \lpsp$ have orienation and will appear in the word $\mathscr{D}(\gamma ')$. If $i\neq i_j$ for any $j$, all simple loops in $V_i$, including ${\gamma_i}'$ are trivial. Therefore ${\gamma_i}'$ has no orienation and will not appear as a letter in $\mathscr{D}(\gamma ')$. If this is the case for all $i$ so that $\mathscr{D}(\beta)$ is the emtpy word, then $\mathscr{D}(\gamma ')$ must also be the empty word $e\in (K')^{-1}(W)$. Suppose on the other hand that $\mathscr{D}(\beta)=\beta_{1}\beta_{2}...\beta_{n}\neq e$ and $i=i_j$ for some $j$. The neighborhood $V_i=V_{i_j}$ was chosen so that all simple loops in $V_{i_j}$ have orientation $\epsilon_j$. Since ${\gamma_i}'\in V_{i_j}$ is simple, it has orienation $\epsilon_j$ and we have ${\gamma_i}'\in V_{i_j}\cap \slpsp^{\epsilon_j}=W_j$. Therefore $\mathscr{D}(\gamma ')={\gamma_{i_1}}'{\gamma_{i_2}}'...{\gamma_{i_n}}'\in W_1W_2...W_n\subseteq (K')^{-1}(W)$.\end{proof}
Using this lemma, we easily prove our main result.\\
\\
\textit{Step 6: Proof of Theorem \ref{mainresult}.}\\
\\
We have already show that $h_{X}$ is a continuous group isomorphism, natural in $X$. To show that $h_{X}$ is also open we suppose $U$ is open in $\fotop$. The commutativity of the diagram in Lemma \ref{diagram} and bijectivity of $h_{X}$ gives equality $$P_{\Omega}^{-1}(h_{X}(U))=\left(R\circ\psi_{X}\circ \rho_{X}\circ P\circ \mathscr{D}\right)^{-1}(U).$$By the previous lemma, the right side of this equation is open in $\lpsp$ and since $P_{\Omega}$ is a quotient map, $h_{X}(U)$ is open in $\pitop$. Therefore, $h_X$ is an open map.
\subsection{The weak suspension spaces $w\sus$ and $\pi_{1}^{top}(w\sus)$}
We pause here to note a deficiency of the suspension spaces $\sus$. It is easy to see that $X$ is compact if and only if there is a countable neighborhood base at the basepoint $x_0$ consisting of neighborhoods of the form $B_n=X\wedge \left[0,\frac{1}{n}\right)\sqcup \left(\frac{n-1}{n},1\right]$. Consequently, if $X$ is a non-compact, first-countable (resp. metric) space, then $\sus$ may not be first countable (resp. a metric space). For this reason we consider a slightly weaker topology on the underlying set of $\sus$, and denote the resulting space as $w\sus$. A basis for the topology of $w\sus$ is given by subsets of the form $V\wedge (a,b)$ and $B_n$, where $V\in \mathscr{B}_{X}$ and $n\geq 2$. Clearly, the identity function $id:\sus\ra w\sus$ is continuous.\\
\\
The "weak suspension" $w\sus$ has a few advantages over $\sus$ including the fact that if $X$ is a subspace of $\mathbb{R}^{n}$, then $w\sus$ may be embedded as a subspace of $\mathbb{R}^{n+1}$.
\begin{example}
For $a\in [0,\infty)$, let$$C_{a}=\{(x,y)\in \mathbb{R}^{2}|(x-a)^{2}+y^{2}=\left(1+a\right)^{2}\}.$$If $X\subset[0,\infty)$ is any subspace of the non-negative real line, then $w \sus$ is homeomorphic to $\bigcup_{a\in X}C_{a}\subset \mathbb{R}^{2}$ with basepoint $(-1,0)$. It is not necessarily true that $\sus\cong \bigcup_{a\in X}C_{a}$ if $X$ is non-compact.
\end{example}
The arguments in this section may be repeated to compute $\pi_{1}^{top}(w \sus)$ or one may prove the following theorem.
\begin{theorem}
The identity map $id:\sus\ra w\sus$, induces a natural isomorphism of quasitopological groups $id:\pitop\ra \pi_{1}^{top}(w\sus)$.
\end{theorem}
\section{Free Topological Groups and $\pitop$}

We now may study the topological properties of $\pitop$ by studying those of $\frpi$ and $\ftpi$. The topological structure of the free topological group $\fty$ on a space $Y$ is often highly complicated and a great deal of effort has gone into determining when a simple description is possible. In particular, the reduction topology of $\fry$ has a remarkably simple description and agrees with that of $\fty$ if and only if $\fry$ is a topological group. We are then confronted with a rather challenging classification problem.\\
\\
\textbf{Question 1:} For which spaces $Y$ is $\fry$ a topological group?\\

According to Proposition \ref{classify1}, $\frpi$ is a topological group whenever $\pitop$ is a topological group. Therefore, a complete answer to Question 1 is necessary (but not sufficient) to characterize the spaces $X$ for which $\pitop$ is a topological group. Question 1 has recieved a good deal of attention for the case when $Y$ is a Tychonoff space. We provide an answer for non-Tychonoff spaces and reference some of the known partial solutions for the Tychonoff case in section 3.2. In doing so, we produce spaces having discontinuous multiplication in their topological fundamental group.
\subsection{The topology of $\fry$ for $Y$ Hausdorff}

Let $Y$ be a Hausdorff space. We study the topological properties of $\fry$ with the goal of characterizing the topological properties of $\pitop$. 
%
\begin{definition} Suppose $w=y_{1}^{\epsilon_1}\dots y_{n}^{\epsilon_n}$ is a non-empty word in $\jy$ and $U_1,\dots,U_n$ is a sequence of open neighborhods in $X$ such that $y_i\in U_i$, $U_j=U_k$ whenever $y_j=y_k$, and $U_j\cap U_k=\emptyset$ whenever $y_j\neq y_k$. An open neighborhood of $w$ in $\jy$ of the form $U=U_{1}^{\epsilon_1}\dots U_{n}^{\epsilon_n}$ is called a \textbf{separating neighborhood of} $\mathbf{w}$. The singleton $\{e\}$ is the only separating neighborhood of the empty word.\end{definition}

Since $Y$ is Hausdorff, every word in $\jy$ has a neighborhood base of separating neighborhoods. We note a subtle distinction in terminology: If $U$ is a separating neighborhood of $w$ and $v\in U$, then $U$ is a neighborhood of $v$ but is
not necessarily a separating neighborhood of $v$.
\begin{lemma} \label{reductionlemma} If $w=y_{1}^{\epsilon_1}\dots y_{n}^{\epsilon_n}$ is a non-empty word in $\jy$, and $U=U_{1}^{\epsilon_1}\dots U_{n}^{\epsilon_n}$ is a separating neighborhood of $w$, then:
\begin{enumerate}
\item $w$ is reduced (i.e. $R(w)=w$) if and only if $U$ contains only reduced words.
\item If $v\in U$, then $|R(w)|\leq |R(v)|\leq |w|$.
\end{enumerate}
\end{lemma}
\begin{proof} 1. follows from the fact that if $y_i\neq y_{i+1}$ or $\epsilon_i=\epsilon_{i+1}$, then $e\notin R\left(U_{i}^{\epsilon_i} U_{i+1}^{\epsilon_{i+1}}\right)$. So if no pairs of consecutive letters can be reduced in $w$, then no pairs of consecutive letters can be reduced in any word in $U$. For 2. if $v\in U$ is reduced, then $|R(w)|\leq|w|=|v|=|R(v)|$ and we are done. One the other hand, if $v$ is not reduced, neither is $w$. Moreover if  $v'=z_{1}^{\epsilon_1}\dots z_{k-1}^{\epsilon_{k-1}}z_{k+2}^{\epsilon_{k+2}}\dots z_{n}^{\epsilon_n}$ is a reduction of $v$ by removing the k-th and (k+1)-th letters of $v$, then there is a corresponding reduction $w'=y_{1}^{\epsilon_1}\dots y_{k-1}^{\epsilon_{k-1}}y_{k+2}^{\epsilon_{k+2}}\dots y_{n}^{\epsilon_n}$ of $w$. The product $U'=U_{1}^{\epsilon_1}\dots U_{k-1}^{\epsilon_{k-1}}U_{k+2}^{\epsilon_{k+2}}\dots U_{n}^{\epsilon_n}$ is then a separating neighborhood of $w'$ containing $v'$. The same arguments apply to $w'$ and $U'$ when $v'$ is reduced or unreduced. Since this will continue until we reach the reduced word of $v$, it follows that $|R(w)|\leq |R(v)|$.\end{proof}

For each integer $n\geq 0$, let $F_{n}(Y)=\{w\in \fry||w|\leq n\}$ denote the set of reduced words of length no greater than $n$. 
\begin{corollary} \label{fnisclosed}
If $Y$ is Hausdorff and $n\geq 0$, then $F_{n}(Y)$ is closed in $\fry$.
\end{corollary}
\begin{proof} Suppose $w\in \jy$ such that $|R(w)|> n$. Now take any separating neighborhood $U$ of $w$. The previous lemma asserts that if $v\in U$, then $|R(v)|\geq |R(w)|> n$. Consequently, $w\in U\subset \jy-R^{-1}(F_{n}(Y))$ proving that $R^{-1}(F_{n}(Y))$ is closed in $\jy$. Since $R$ is a quotient map $F_{n}(Y)$ is closed in $\fry$.\end{proof}

We now find some properties of $\fry$ which are often desirable in free topological groups. Let $Z$ denote the set of all finite sequences $\zeta=\epsilon_1,...,\epsilon_n$ with $\epsilon_i\in \{\pm 1\}$, including the empty sequence. For each $\zeta=\epsilon_1,...,\epsilon_n\in Z$, let $Y^{\zeta}=\{y_{1}^{\epsilon_1}\dots y_{n}^{\epsilon_n}|y_i\in Y\}$ and recall that $\jy=\coprod_{\zeta\in Z}Y^{\zeta}$. Let $|\zeta|$ denote the length of each sequence and $R_{\zeta}:Y^{\zeta}\ra F_{|\zeta|}(Y)$ be the restriction of the reduction quotient map. The next two propositions correspond exactly to those in Statement 5.1 of \cite{Sipacheva}.
\begin{proposition} \label{inductivelimit} If $Y$ is Hausdorff, $\fry$ has the inductive limit topology of the sequence of closed subspaces $\{F_{n}(Y)\}_{n\geq 0}$, i.e. $\fry\cong \varinjlim_{n}F_{n}(Y)$.
\end{proposition}
\begin{proof}  Suppose $C\subseteq \fry$ such that $C\cap F_{n}(Y)$ is closed in $F_{n}(Y)$ for each $n\geq 0$. But $R^{-1}(C)\cap Y^{\zeta}=R_{\zeta}^{-1}(C)=R_{\zeta}^{-1}(C\cap F_{|\zeta|}(Y))$ is closed in $Y^{\zeta}$ for each $\zeta$, since $R_{\zeta}$ is continuous. Since $\jy$ is the disjoint union of the $Y^{\zeta}$, $R^{-1}(C)$ is closed in $\jy$ and $R$ being quotient gives that $C$ is closed in $\fry$.\end{proof}
\begin{proposition} \label{restrictedquotient} If $Y$ is Hausdorff and $n\geq 0$, the restriction $R_n:\coprod_{i=0}^{n}(Y\sqcup Y^{-1})^{i}\ra F_{n}(Y)$ of the reduction homomorphism is a quotient map.
\end{proposition}
\begin{proof} Suppose $A\subseteq F_{n}(Y)$ such that $R_{n}^{-1}(A)$ is closed in $\coprod_{i=0}^{n}(Y\sqcup Y^{-1})^{i}=\coprod_{|\zeta|\leq n}Y^{\zeta}$. Since $F_{n}(Y)$ is closed in $\fry$ and $R$ is a quotient map, it suffices to show that $R_{\zeta}^{-1}(A)=\{a=a_{1}^{\epsilon_1}...a_{k}^{\epsilon_{k}}|R(a)\in A\}$ is closed in $Y^{\zeta}$ for each $\zeta=\epsilon_1,...,\epsilon_k\in Z$. If $|\zeta|\leq n$, then $R_{\zeta}^{-1}(A)=R_{n}^{-1}(A)\cap Y^{\zeta}$ is closed by assumption. For $|\zeta|>n$, we proceed by induction, and suppose $R_{\delta}^{-1}(A)$ is closed in $Y^{\delta}$ for all $\delta$ of length $|\delta|=n,n+1,...,|\zeta|-1$. Let $w=y_{1}^{\epsilon_1}...y_{k}^{\epsilon_{k}}\in Y^{\zeta}$ such that $R(w)\notin A$. If $w$ is reduced, $U$ is a separating neighborhood of $w$, and $v\in U$, then by Lemma \ref{reductionlemma}, $n<|\zeta|=|R(w)|=|R(v)|$. Thus $R(U)\cap F_{n}(Y)=\emptyset$ and since $A\subseteq F_{n}(Y)$ we have $U\cap R_{\zeta}^{-1}(A)=\emptyset$. Therefore we may suppose that $w$ is not reduced. For each $i\in \{1,...,k-1\}$ such that $y_{i}=y_{i+1}$ and $\epsilon_i=-\epsilon_{i+1}$ we find a separating neighborhood $U_i$ of $w$ in the following way. Remove the two letters $y_{i}^{\epsilon_i},y_{i+1}^{\epsilon_{i+1}}$ from $w$ to obtain a word $w_{i}=y_{1}^{\epsilon_1}...y_{i-1}^{\epsilon_{i-1}} y_{i+2}^{\epsilon_{i+2}}...y_{k}^{\epsilon_{k}}$ such that $R(w_{i})=R(w)\notin A$. Let $\zeta_{i} =\epsilon_1,...,\epsilon_{i-1},\epsilon_{i+2},...,\epsilon_k$ so that $|\zeta_{i} |=|\zeta|-2$ and $w_{i}\in Y^{\zeta_{i} }-R_{\zeta_{i} }^{-1}(A)$. By our induction hypothesis $R_{\zeta_{i} }^{-1}(A)$ is closed in $Y^{\zeta_{i} }$ and so we may find a separaing neighborhood $V_i=A_{1}^{\epsilon_1}...A_{i-1}^{\epsilon_{i-1}} A_{i+2}^{\epsilon_{i+2}}...A_{k}^{\epsilon_{k}}$ of $w_i$ such that $V_i\cap R_{\zeta_{i} }^{-1}(A)=\emptyset$. We may then find a neighborhood $A_{i}=A_{i+1}$ of $y_i=y_{i+1}$ such that $$U_i=A_{1}^{\epsilon_1}...A_{i-1}^{\epsilon_{i-1}}A_{i}^{\epsilon_i} A_{i+1}^{\epsilon_{i+1}}A_{i+2}^{\epsilon_{i+2}}...A_{k}^{\epsilon_{k}}$$is a separating neighborhood of $w$. Now take a separating neighborhood $U$ of $w$ such that $w\in U\subseteq \bigcap_{i}U_i$ (the intersection ranges over $i$ such that $y_{i}=y_{i+1}$ and $\epsilon_i=-\epsilon_{i+1}$). It now suffices to show that $R(v)\notin A$ whenever $v=z_{1}^{\epsilon_1}...z_{k}^{\epsilon_k}\in U$. If $v$ is reduced, then $n<|\zeta|=|w|=|R(v)|$. Thus $R(v)\notin F_{n}(Y)$ and we are done. If $v$ is not reduced, then there is an $i_0\in \{1,...,k-1\}$ such that $z_{i_0}=z_{{i_0}+1}$ and $\epsilon_{i_0}=-\epsilon_{{i_0}+1}$. If $v_{i_0}=z_{1}^{\epsilon_1}...z_{{i_0}-1}^{\epsilon_{{i_0}-1}} z_{{i_0}+2}^{\epsilon_{{i_0}+2}}...z_{k}^{\epsilon_{k}}$ is the word obtained by removing $z_{i_0}^{\epsilon_{i_0}},z_{{i_0}+1}^{\epsilon_{{i_0}+1}}$, we have $R(v)=R(v_{i_0})$. But since $U$ is a separating neighborhood of $w$ we must also have $y_{i_0}=y_{{i_0}+1}$. Since $U\subseteq U_{i_0}$, we have $v_{i_0}\in V_{i_0}$ and therefore $R(v)=R(v_{i_0})\notin A$.\end{proof}
\begin{theorem} \label{propsoffry}
The following are equivalent:
\begin{enumerate}
\item $Y$ is Hausdorff.
\item $\fry$ is $T_1$.
\item For each $n\geq 1$, the canonical map $\sigma_{n}:Y^{n}\ra \fry$ taking $(y_1,...,y_n)$ to the word $y_1...,y_n$ is a closed embedding.
\end{enumerate}
\end{theorem}
\begin{proof} 1. $\Rightarrow$ 2. If $Y$ is Hausdorff, then the singleton $F_{0}(Y)=\{e\}$ containing the identity is closed. Since $\fry$ is a quasitopological group, it is $T_1$.\\
\\
2. $\Rightarrow$ 1. If $Y$ is not Hausdorff, then we suppose $y$ and $z$ are distinct points in $Y$ which cannot be separated by disjoint open sets. If $W$ is any open neighborhood of reduced word $yz^{-1}$ in $\fry$, then there is an open neighborhood $U$ of $y$ and $V$ of $z$ in $Y$ such that $yz^{-1}\in UV^{-1}\subset R^{-1}(W)\subset \jy$. But there is a $w\in U\cap V$ by assumption and so $UV^{-1}$ contains the word $ww^{-1}$. Therefore $e=R(ww^{-1})\in W$. But if every neighborhood of $yz^{-1}$ in $\fry$ contains the identity, $\fry$ is not $T_1$.\\
\\
1. $\Rightarrow$ 3. Suppose $A$ is a closed subspace of $Y^n$. Let $j:Y^n\hookrightarrow \coprod_{i=0}^{n}(Y\sqcup Y^{-1})^{i}$ be the map $j(y_1,...,y_n)=y_1...y_n$ so that $R_n\circ j=\sigma_{n}$. Since $j$ is a closed embedding, $R_{n}^{-1}(\sigma_{n}(A))=j(A)$ is closed in $\coprod_{i=0}^{n}(Y\sqcup Y^{-1})^{i}$. But $R_n$ is quotient and $F_{n}(Y)$ is closed in $\fry$. Therefore $\sigma_{n}(A)$ is closed in $\fry$.\\
\\
3. $\Rightarrow$ 1. If $Y$ is not Hausdorff, the argument for 2. $\Rightarrow$ 1. asserts that there are distinct $y,z\in Y$ such that any open neighborhood in $\jy$ of the three letter word $zyz^{-1}$ contains a word $zww^{-1}$ which satisfies $R(zww^{-1})=z$. Therefore, any open neighborhood of reduced word $zyz^{-1}$ in $\fry$ contains the one letter word $z\in \sigma_{1}(Y)$. This means the image of $\sigma_{1}$ is not closed in $\fry$ and $\sigma_{1}$ cannot be a closed embedding.\end{proof}
Since we are primarily interested in T1 quasitopological groups, this theorem motivates our study of $\fry$ for Hausdorff spaces $Y$. It also provides an immediate answer to Question 1 for Hausdorff, non-completely regular spaces.
\begin{corollary} \label{notatopologicalgroup}
If $Y$ is Hausdorff but not completely regular, then $\fry$ is not a topological group.
\end{corollary}
\begin{proof} If $Y$ is Hausdorff and $\fry$ is a topological group, then $\fry\cong \fty$ and $\sigma:Y\hookrightarrow \fty$ is an embedding. In \cite{Thomas} it is shown that $\sigma:Y\ra \fty$ is an embedding if and only if $Y$ is completely regular. Therefore, $Y$ must also be completely regular.\end{proof}
Recall that $\frpi$ is a topological group, whenever $\pitop$ is a topological group. Therefore, if $\piztop$ is Hausdorff but not completely regular, then $\pitop$ is not a topological group. This fact enables us to produce explicit counterexamples to the assertion that topological fundamental groups are always a topological group.
\begin{example}
Let $\mathbb{Q}_{K}$ denote the rational numbers with the subspace topology of the real line with the K-topology. $\pi_{0}^{top}(\mathbb{Q}_{K})\cong \mathbb{Q}_{K}$ is totally path disconnected, Hausdorff but is not regular. Therefore $\pi_{1}^{top}(\Sigma (\mathbb{Q}_{K})_{+})\cong F_{R}(\mathbb{Q}_{K})$ is not a topological group.
\end{example}
\begin{example}
For every Hausdorff, non-completely regular space $Y$, there is, by Remark \ref{essentialsurj}, a paracompact Hausdorff space $X=\mathcal{H}(Y)$ such that $\piztop\cong Y$. Since $\fry\cong \frpi$ is not a topological group, then neither is $\pitop$.
\end{example}

We may also apply Theorem \ref{propsoffry} to fully characterize when all powers of a continuous surjection are quotient.
\begin{theorem} \label{quotientmappowers}
Let $X$ and $Y$ be spaces with $X$ Hausdorff and $q:X\ra Y$ be a continuous surjection. The induced continuous, epimorphism $F_{R}(q):\frx\ra \fry$ is a quotient map if and only if the powers $q^n:X^n\ra Y^n$ are quotient maps for all $n\geq 1$.
\end{theorem}
\begin{proof} If $q^n:X^n\ra Y^n$ is a quotient map for each $n\geq 1$, then so is $M_{T}^{\ast}(q)=\coprod_{n\geq 0}(q\sqcup q)^{n}:\jx\ra \jy$. Since the diagram $$\xymatrix{ \jx \ar[r]^{M_{T}^{\ast}(q)} \ar[d]_{R_X} & \jy \ar[d]^{R_Y} \\
\frx \ar[r]_{F_{R}(q)} & \fry }$$commutes and the reduction maps are quotient maps, $F_{R}(q)$ is a quotient map.\\
\\
To prove the converse, we let $\sigma^{X}_{n}:X^n\ra \frx$ and $\sigma^{Y}_{n}:Y^n\ra \fry$ be the canonical maps of Theorem \ref{propsoffry} and $\tilde{X}^{n}=\sigma^{X}_{n}(X^n)$ and $\tilde{Y}^{n}=\sigma^{Y}_{n}(Y^n)$ be their images. Since $X$ is Hausdorff, $\sigma^{X}_{n}$ is a closed embedding. We show the restriction $p=F_{R}(q)|_{\tilde{X}^{n}}:\tilde{X}^{n}\ra \tilde{Y}^{n}$ is a quotient map using the commutative diagram$$\xymatrix{ X^n \ar[d]_{q^n} \ar[r]^{\sigma_{X}^{n}} & \frx \ar[d]^{F_{R}(q)}  & \jx \ar[d]^{M_{T}^{\ast}(q)} \ar[l]_{R_{X}} \\ Y^n \ar[r]_{\sigma_{Y}^{n}} & \fry  & \jy \ar[l]^{R_{Y}} }$$where the reduction maps are distinguished with subscripts. To see that $p$ being quotient implies $q^n$ is quotient, take $C\subseteq Y^n$ such that $(q^n)^{-1}(C)$ is closed in $X^n$. Then $\sigma_{n}^{X}((q^n)^{-1}(C))=p^{-1}(\sigma_{n}^{Y}(C))$ is closed in $\tilde{X}^{n}$ and consequently $\sigma_{n}^{Y}(C)$ is closed in $\tilde{Y}^{n}$. Since $\sigma_{n}^{Y}$ is a continuous injection, $C$ is closed in $Y^n$.\\
\\
Suppose $A\subseteq \tilde{Y}^{n}$ such that $p^{-1}(A)$ is closed in $\tilde{X}^{n}$. Since $R_X$ and $F_{R}(q)$ are assumed to be quotient and $\tilde{Y}^{n}$ is closed in $\fry$, it suffices to show that $$B^{\zeta}=R_{X}^{-1}(F_{R}(q)^{-1}(A)\cap X^{\zeta}=\{x=x_{1}^{\epsilon_1}...x_{k}^{\epsilon_k}|R_{Y}(M_{T}^{\ast}(q)(x))= R_{Y}(q(x_1)^{\epsilon_1}...q(x_k)^{\epsilon_k})\in A\}$$ is closed in $X^{\zeta}$ for each $\zeta=\epsilon_1,...,\epsilon_k$. We proceed by induction on $|\zeta|=k$. It is clear that if $|\zeta|<n$, then $B^{\zeta}=\emptyset$. Additionally, if $|\zeta|=n$ and $\zeta\neq 1,1,...,1$, then $B^{\zeta}=\emptyset$. On the other hand, if $|\zeta|=n$ and $\zeta= 1,1,...,1$, then $B^{\zeta}=\{x_1...x_n|q(x_1)...q(x_n)\in A\}=R_{X}^{-1}(p^{-1}(A))\cap X^{\zeta}$ is closed by assumption. Now we suppose that $|\zeta|>n$ and $B^{\zeta}$ is closed in $X^{\delta}$ for all $\delta$ such that $|\delta|=n,n+1,...,|\zeta|-1$. Let $x=x_{1}^{\epsilon_1}...x_{k}^{\epsilon_k}\in X^{\zeta}-B^{\zeta}$ and $y=\mathcal{J}^{\pm}(q)(x)=q(x_{1})^{\epsilon_1}...q(x_{k})^{\epsilon_k}$. Since $x\notin B^{\zeta}$, we have $R_{Y}(y)\notin A$. Let $E=E_{1}^{\epsilon_1}...E_{k}^{\epsilon_k}$ be a separating neighborhood of $y$ in $\jy$. Since $M_{T}^{\ast}(q)$ is continuous, there is a separating neighborhood $D=D_{1}^{\epsilon_1}...D_{k}^{\epsilon_k}$ of $x$, such that $q(D_i)\subseteq E_i$ for each $i\in \{1,...,k\}$. Since $E$ is a separating neighborhood, if $q(x_i)\neq q(x_{j})$, then $q(D_i)\cap q(D_j)=\emptyset$. Now we consider the cases when $y$ is and is not reduced.\\
\\
If $y$ is reduced and $v\in D$, then $M_{T}^{\ast}(q)(v)\in E$ must also be reduced by Lemma \ref{reductionlemma}. Therefore $n<|\zeta|=|y|= |R_{Y}(M_{T}^{\ast}(v))|$, i.e. the reduced word of $M_{T}^{\ast}(q)(v)$ has length greater than $n$ and cannot lie in $A\subseteq \tilde{Y}^{n}$. Therefore $D\cap B^{\zeta}=\emptyset$.\\
\\
If $y$ is not reduced, then for each $i\in \{1,...,k-1\}$ such that $q(x_{i})=q(x_{i+1})$ and $\epsilon_{i}=-\epsilon_{i+1}$, we let $w_i=q(x_{1})^{\epsilon_1}\dots q(x_{i-1})^{\epsilon_{i-1}}q(x_{i+2})^{\epsilon_{i+2}}\dots q(x_{k})^{\epsilon_k}\in \jy$ and $u_i=x_{1}^{\epsilon_1}\dots x_{i-1}^{\epsilon_{i-1}}x_{i+2}^{\epsilon_{i+2}}\dots x_{k}^{\epsilon_k}$ be the words obtained by removing the i-th and (i+1)-th letters from $y$ and $x$ respectively. We also let $\zeta_i=\epsilon_1,\dots,\epsilon_{i-1},\epsilon_{i+2},...\epsilon_{k}$. This gives $F_{R}(q)(R_{X}(u_i))=R_{Y}(M_{T}^{\ast}(u_i))=R_{Y}(w_i)=R_{Y}(y)\notin A$ and consequently $u_i\in X^{\zeta_i}-B^{\zeta_i}$. We know by our induction hypothesis that $B^{\zeta_i}$ is closed in $X^{\zeta_i}$ and so we may find a separating neighborhood $V_i= A_{1}^{\epsilon_1}\dots A_{i-1}^{\epsilon_{i-1}}A_{i+2}^{\epsilon_{i+2}}\dots A_{k}^{\epsilon_k}$ of $u_i$ contained in $X^{\zeta_i}-B^{\zeta_i}$. Let $A_i=A_{i+1}=X$ so that $$U_i=A_{1}^{\epsilon_1}...A_{i-1}^{\epsilon_{i-1}} A_{i}^{\epsilon_i}A_{i+1}^{\epsilon_{i+1}}A_{i+2}^{\epsilon_{i+2}}...A_{k}^{\epsilon_k}$$ is an open neighborhood of $x$. Now take a separating neighborhood $U$ of $x$ such that $U\subseteq D\cap\bigcap_{i}U_i$ where the intersection ranges over the $i\in\{1,\dots ,k-1\}$ such that $q(x_{i})=q(x_{i+1})$ and $\epsilon_{i}=-\epsilon_{i+1}$. It now suffices to show that $F_{R}(q)(R_X(v))=R_{Y}(q(z_{1})^{\epsilon_1}\dots q(z_{k})^{\epsilon_k})\notin A$ whenever $v=z_{1}^{\epsilon_1}\dots z_{k}^{\epsilon_k}\in U$. If $M_{T}^{\ast}(q)(v)=q(z_{1})^{\epsilon_1}\dots q(z_{k})^{\epsilon_k}$ is reduced, then $n<|\zeta|=|x|=|R_{Y}(\mathcal{J}^{\pm}(q))|$ and $R_{Y}(M_{T}^{\ast}(q)(v))\notin A$. On the other hand, suppose $q(z_{1})^{\epsilon_1}\dots q(z_{k})^{\epsilon_k}$ is not reduced. There is an $i_0\in\{1,...,k-1\}$ such that $q(z_{i_0})=q(z_{i_{0}+1})$ and $\epsilon_{i_0}=-\epsilon_{i_{0}+1}$. But $z_{i_0}\in D_{i_0}$ and $z_{{i_0}+1}\in D_{{i_0}+1}$, so we must have $q(x_{i_0})=q(x_{{i_0}+1})$. Since $v\in U\subseteq U_{i_0}$, we have $$v_{i_0}=z_{1}^{\epsilon_1}\dots z_{{i_0}-1}^{\epsilon_{{i_0}-1}}z_{{i_0}+2}^{\epsilon_{{i_0}+2}}\dots z_{k}^{\epsilon_k}\in V_{i_0}\subseteq X^{\zeta_{i_0}}-B^{\zeta_{i_0}}$$Therefore $$F_{R}(q)(R_X(v))=R_{Y}(M_{T}^{\ast}(q)(v))=R_{Y}(M_{T}^{\ast}(q)(v_{i_0}))=F_{R}(q)(R_X(v_{i_0}))\notin A$$proving that $U\cap B^{\zeta}=\emptyset$ and $B^{\zeta}$ is closed in $X^{\zeta}$.\end{proof}
Applying this characterization to our computation of $\pitop$ we see that
\begin{theorem} \label{classification1}
If $X$ is Hausdorff, then $\pitop\cong\fotop\cong \frpi$ if and only if $(P_X)^n:X^n\ra \piztop^n$ is a quotient map for all $n\geq 1$.
\end{theorem}
\begin{proof}
If $(P_X)^n:X^n\ra \piztop^n$ is a quotient map for all $n\geq 1$, then $\pitop\cong\fotop\cong \frpi$ by Theorem \ref{mainresult} and Proposition\ref{classify1}. Conversely, if $\fotop\cong \frpi$, then $F_{R}(P_X)$ is quotient by Proposition \ref{classify1}. By Theorem \ref{quotientmappowers}, $(P_X)^n:X^n\ra \piztop^n$ must be a quotient map for all $n\geq 1$.
\end{proof}

Theorem \ref{propsoffry} also allows us to remark on the countability of $\fry$. The argument used to prove the next statement is based on that used by Fabel \cite{F8} to show that the Hawaiian earring group is not first countable. For a word $w\in \jy$ and letter $y\in Y$, we let $O_{y}(w)$ denote the numbers of times the letter $y$ or $y^{-1}$ appears in $w$. Note that for every word $w\in \jy$ and $y\in Y$, $O_{y}(w)\geq O_{y}(R(w))$.
\begin{theorem} \label{firstcountability}
If $Y$ is Hausdorff, then $Y$ is discrete if and only if $\fry$ is first countable
\end{theorem}
\begin{proof} If $Y$ is discrete, then $\fry$ is discrete and thus first countable. Suppose that $Y$ is a non-discrete Hausdorff space (so $\fry$ is $T_1$) such that $\fry$ is first countable. We let $(J,\leq)$ denote the directed set of pointwise coverings $(Cov(Y),\preceq)$ as in section 1. If $Y$ is not discrete, then $J$ has no maximal element but does have a minimal element $M=\{V^y\}_{y\in Y}$ (the covering where $V^y=Y$ for each $y\in Y$). We also have a convergent net $y_j\shortrightarrow y_0$ in $Y$ where $y_0\neq y_j$ for each $j\in J$. Suppose $\{B_1,B_2,...\}$ is a countable basis of open neighborhoods at the identity $e$ in $\fry$ where $B_{i+1}\subseteq B_{i}$ for each $i$. We note that the reduced word of the n-fold concatenation $w_n=(y_0y_{M}y_{0}y_{0}^{-1}y_{M}^{-1}y_{0}^{-1})^{n}$ is the empty word $e$ for each $n\geq 1$ so $w_n\in R^{-1}(B_i)$ for all $i,n\geq 1$. Since $Y$ is non-discrete Hausdorff, we may find $M<j_1<k_1<j_{2}<k_{2}<j_{3}<...$ in $J$ such that $y_M$ and all of the $y_{j_n}$ and $y_{k_n}$ are distinct. We now consider the reduced word $$v_n=\left(y_{j_n}y_{M}y_{k_n}{y_{j_n}}^{-1}{y_{M}}^{-1}{y_{k_n}}^{-1}\right)^{n}\in R^{-1}(B_n)$$noting that $O_{y_M}(v_n)=2n$ for each $n\geq 1$. The set $C=\{v_n|n \geq 1\}\subset \fry$ does not contain the identity. We show that $C$ is closed in $\fry$ by showing $R^{-1}(C)$ is closed in $\jy$. If this is done, $v_n$ cannot be a sequence converging to the empty word.\\
\\
Suppose that $\{z_{k}\}_{k\in K}\shortrightarrow z$ is a convergent net of words in $\jy$ such that $z_k\in R^{-1}(C)$ but such that $z\notin R^{-1}(C)\cap Y^{\zeta}$ for some $\zeta\in Z$. Since $\jy=\coprod_{\zeta\in Z}Y^{\zeta}$, there is a $k_0\in K$ such that $z_k\in Y^{\zeta}$ (and consequently $|z_{k}|=|z|$) for every $k\geq k_0$. But since $R(z_k)\in C$ for each $k\in K$, we may write $R(z_{k})=v_{n_k}$. Suppose the net of integers $n_k$ is bounded by integer $N$. Then $R(z_k)\in \{v_1,v_2,...,v_N\}$ for each $k\in K$. But since $\fry$ is $T_1$ the finite set $\{v_1,v_2,...,v_N\}$ is closed in $\fry$ and we must have $R(z)\in \{v_1,v_2,...,v_N\}\subseteq C$. This, however, is a contradiction. Suppose, on the other hand, that $n_k$ is unbounded.
Since $R(z_{k})=R( v_{n_{k}})=v_{n_k}$ for each $k\in K$, we have that $$|z_{k}|\geq O_{y_M}\left(z_{k}\right)\geq O_{y_M}\left(v_{n_{k}}\right)=2n_{k}$$for each $k \in K$. But this means $|z_{k}|$ is an unbounded net of integers, contradicting the fact that eventually $|z_{k}|=|z|$. Therefore we must have that $R(z)\in C$ so that $C$ is closed in $\fry$.\end{proof}
Our previous characterization of discreteness and the fact that the topology of $\pitop$ is finer than that of $\frpi$ then gives
\begin{corollary}
If $\piztop$ is Hausdorff, then $\pitop$ is first countable if and only if $\piztop$ is discrete.
\end{corollary}
\subsection{Separation Properties in $\pitop$}
We note some immediate consequences of the previous sections. Specifically, the fact that $u_{\ast}:\piztop\ra \pitop$ is a continuous injection implies that
\begin{proposition} \label{hausdorff}
If $\pitop$ is T1 (resp. Hausdorff) then so is $\piztop$.
\end{proposition}
Moreover, if $\piztop$ is Hausdorff, then $\sigma:\piztop\ra \frpi$ is a closed embedding. Since the topology of $\pitop$ is finer than that of $\frpi$, the map $u_{\ast}:\piztop\ra \pitop$ is also a closed embedding. When then have
\begin{proposition} \label{properties} If $\piztop$ is Hausdorff and (Pr) is a topological property hereditary to closed subspaces, then $\piztop$ has property (Pr) whenever $\pitop$ has property (Pr).
\end{proposition}
It is natural to ask for which properties (Pr) the converse also holds. For instance, while it is clear that $\pitop$ Hausdorff $\Rightarrow$ $\piztop$ Hausdorff, the converse is not. To see this, it suffices to construct a Hausdorff space $Y$ which fails to have the following property $(H')$.\\
\\
$(H')$: For every pair of distinct points $a,b\in Y$, there are open neighborhoods $U$ of $a$ and $V$ of $b$ such that for every $y\in Y$, there is a neighborhood $W$ of $y$ such that $W$ does not intersect both $U$ and $V$.\\

It can be shown quite easily that $$Y\text{ is Regular and }T_1\text{ }\Rightarrow\text{ }Y\text{ has property }(H')\text{ }\Rightarrow\text{ }Y\text{ is Hausdorff}$$
\begin{proposition}
If $\fry$ is Hausdorff, then $Y$ must have property $(H')$.
\end{proposition}
\begin{proof} Suppose $a,b\in Y$ are two distinct points failing to satisfy the condition in $(H')$. Suppose $N$ and $M$ are open neighborhoods of  $ab^{-1}$ and $e$ respectively in $\fry$. Since $\fry$ has the quotient topology of $\jy$, there are open neighborhoods $U$ of $a$ and $V$ of $B$ in $Y$ such that $ab^{-1}\in UV^{-1}\subseteq R^{-1}(N)$. By assumption, there is a $y\in Y$ such that for every neighborhood $W$ of $y$, $W\cap U\neq \emptyset$ and $W\cap V\neq\emptyset$. But $R(yy^{-1})=e$ and so there is an open neighborhood $W$ of $y$ such that $yy^{-1}\in WW^{-1}\subseteq R^{-1}(M)$. But there is a $c\in W\cap U$ and $d \in W\cap V$ so $cd^{-1}\in UV^{-1}$ and $cd^{-1}\in WW^{-1}$. Therefore $R(cd^{-1})\in M\cap N$ and $\fry$ cannot be Hausdorff.\end{proof}
This indicates a potential characterization of spaces for which $\fry$ is Hausdorff.\\
\\
\textbf{Question 2:} Is $\fry$ Hausdorff whenever $Y$ has property $(H')$?\\

We now give an explicit example of a totally path disconnected, Hausdorff space $X$ which does not have property $(H')$. In this case $\pitop\cong\frx$ is T1 but fails to be Hausdorff (and consequently to be a topological group).
\begin{example} \label{nonhausdorffpi1top} \emph{
We begin by defining the underlying set of a space $X$. Let $K=\{\frac{1}{n}|n\geq 1\}$, $-K=\{-k|k\in K\}$, and $X=(K\times (-K\sqcup \{0\}\sqcup K))\sqcup \{a,b\}$. We define a basis for the topology of $X$ as follows. For each $(r,s)\in K\times (-K\sqcup K)$, the singleton $\{(r,s)\}$ is open in $X$. Let $K_m$ be the set $K_m=\{\frac{1}{k}|k\geq m\}$ for each integer $m\geq 1$. A basic open neighborhood of $(r,0)$ is of the form $\{r\} \times (-K_{m}\sqcup \{0\}\sqcup K_m)$, of $a$ is of the form $\{a\}\sqcup (K_m\times K)$, of $b$ is of the form $\{b\}\sqcup (K_m\times (-K))$. This space is Hausdorff and totally path disconnected but any pair of neighborhoods of $a$ and $b$ are both intersected by any neighborhood of $(r,0)$ for some $r\in K$.}
\end{example}
For the property of being functionally Hausdorff, the converse of Proposition \ref{properties} does hold. Recall that a space $Y$ is functionally Hausdorff if for each pair of distinct points $a,b\in Y$, there is a real valued function $f:Y\ra \mathbb{R}$ such that $f(a)\neq f(b)$.
\begin{proposition} $\pitop$ is functionally Hausdorff if and only if $\piztop$ is functionally Hausdorff.
\end{proposition}
\begin{proof} The fact that $\piztop$ always injects continuously into $\pitop$ indicates that $\piztop$ is functionally Hausdorff whenever $\pitop$ is. Conversely, suppose $\piztop$ is functionally Hausdorff, Theorem 2.3 of \cite{Thomas} asserts that $\piztop$ is functionally Hasudorff if and only if $\ftpi$ is functionally Hausdorff. Since the topology of $\pitop$ is finer than the topology of $\ftpi$, $\pitop$ is functionally Hausdorff whenever $\ftpi$ is.\end{proof}
\subsection{When is $\pitop$ a topological group?}
Corollary \ref{notatopologicalgroup} of Section 3.1 reduces Question 1 to:\\
\\
\textbf{Question 3:} For which Tychonoff spaces $Y$ is $\fry$ a topological group?\\
\\
This question has close ties to many open classification problems in the study of free topological groups \cite{Sipacheva}. If we compose the continuous identity homomorphism $\phi_{Y}:\fry\ra \fty$ with reduction we obtain a map which may still be thought of as reduction but which is commonly denoted as$$\xymatrix{\mathbf{i}:\jy \ar[r]^{R} & \fry \ar[r]^{\phi_{Y}} & \fty }.$$Typically, the restriction of $\mathbf{i}$ to the subspace $\coprod_{i=0}^{n}(Y\sqcup Y^{-1})^{i}$ is denoted as $\mathbf{i}_{n}:\coprod_{i=0}^{n}(Y\sqcup Y^{-1})^{i}\ra F_{n}(Y)$. It is a well known fact that the subsets $F_{n}(Y)$ are closed in $\fty$ for a Tychonoff space $Y$. Part of the following proposition appears as Statement 5.1 of \cite{Sipacheva}.
\begin{theorem} \label{tychonoffconditions}
For a Hausdorff space $Y$, the following are equivalent:
\begin{enumerate}
\item $\fry$ is a topological group.
\item $\mathbf{i}$ is a quotient map.
\item 3a. For each $n\geq 1$, $\mathbf{i}_{n}$ is a quotient map and 3b. $\fty$ is the inductive limit of the subspaces $F_{n}(Y)$ (i.e. $\fty\cong \varinjlim_{n}F_{n}(Y)$).
\end{enumerate}
\end{theorem}
\begin{proof} 1. $\Leftrightarrow$ 2. is obvious since $R$ is a quotient map. 1. $\Rightarrow$ 3. follows easily from Propositions \ref{inductivelimit} and \ref{restrictedquotient}. To show 3. $\Rightarrow$ 2., we suppose $A\subseteq \fty$ such that $\mathbf{i}^{-1}(A)$ is closed in $\jy$. For each $n$ we have that $$\mathbf{i}_{n}^{-1}(A\cap F_{n}(Y))=\mathbf{i}^{-1}(A)\cap \coprod_{i=0}^{n}(Y\sqcup Y^{-1})^{i}=\coprod_{|\zeta|\leq n}(\mathbf{i}^{-1}(A)\cap Y^{\zeta})$$ is closed since it is a finite, disjoint union of closed subspaces. Since $\mathbf{i}_n$ is assumed to be quotient for each $n$, $A\cap F_{n}(Y)$ is closed in $F_{n}(Y)$ for each $n$. Since $\fty$ is inductive limit of the $F_{n}(Y)$, $A$ is closed in $\fry$.\end{proof}
Of course, Corollary \ref{notatopologicalgroup} implies that if $Y$ is not completely regular, then all three conditions of the previous theorem fail. This justifies the effort to characterize only Tychonoff spaces for which 3a. and 3b. are individually satisfied. A full characterization remains elusive on both counts and a brief survey of the known partial characterizations may be found in \cite{Sipacheva}.\\

Using Theorems \ref{classification1} and \ref{tychonoffconditions}, we summarize our characterization of spaces for which $\pitop$ is a Hausdorff topological group in the following theorem. It boils down to a simple condition on the path component space and three well known conditions in topology.
\begin{theorem} \label{classificationtheorem}
Suppose $X$ is Hausdorff. Then $\pitop$ is a Hausdorff topological group if and only if the following conditions hold:
\begin{enumerate}
\item $\piztop$ is Hausdorff.
\item $(P_X)^{n}$ is a quotient map for every $n\geq 1$.
\item $\mathbf{i}_{n}:\coprod_{i=0}^{n}(\piztop\sqcup \piztop^{-1})^{i}\ra F_{n}(\piztop)$ is a quotient map for every $n\geq 1$.
\item $\ftpi$ has the inductive limit topology of the subspaces $F_{n}(\piztop)$.
\end{enumerate}
Moreover, if $\piztop$ is Hausdorff but fails to be completely regular, then conditions 2.-4. fail.
\end{theorem} 
\begin{proof}
If $\pitop$ is Hausdorff, then by Proposition \ref{hausdorff}, $\piztop$ is Hausdorff. If $\pitop$ is also a topological group, then $\pitop\cong \fotop\cong\frpi\cong\ftpi$. The second topological isomorphism and Theorem \ref{classification1} give that $(P_X)^{n}$ is a quotient map for every $n\geq 1$. The last topological isomorphism and Theorem \ref{tychonoffconditions} give that conditions 3. and 4. hold. Conversely, suppose conditions 1.-4. hold. But 1.,3.,4. imply $\frpi\cong\ftpi$ by Theorem \ref{tychonoffconditions}. Additionally, 2. implies $\fotop\cong \frpi$ by Corollary \ref{classify1}. Therefore $\pitop\cong \fotop\cong \frpi\cong \ftpi$ is a topological group. The last statement is direct consequence of Corollary \ref{notatopologicalgroup}.
\end{proof}
Conditions 3.-4. have recieved a good deal of interest. Regarding 2., there are many known characterizations of when products of quotient maps are again quotient maps. For instance, it is known that the functor $X\times -$ preserves quotients if and only if $X$ is core-compact. The author is unaware of a general characterization of quotient maps $q:X\ra Y$ such that all powers $q^n:X^n\ra Y^n$ are quotient.\\

In \cite{Graev,MMO} it is shown that if $X$ is a compact Hausdorff space or more generally a Tychonoff, $k_{\omega}$-space (a direct limit of an increasing sequence of compact subspaces), then $\mathbf{i}$ is quotient and $\frx$ is a topological group. The next statement follows directly from this fact and Corollaries \ref{tgclassification2} and \ref{totallypathdisconnected}.
\begin{corollary}
If $X$ is a Tychonoff $k_{\omega}$-space, then $\pitop$ is a topological group. If, in addition, $X$ is totally path disconnected, then $\pitop\cong \ftx$.
\end{corollary}
\begin{example} \emph{
Let $\mathbb{N}\cup\{\infty\}$ be the one-point compactification of the discrete space of natural numbers. Since $\mathbb{N}\cup\{\infty\}$ is totally path disconnected, $\pi_{1}^{top}(\Sigma((\mathbb{N}\cup\{\infty\})_{+}))\cong F_{R}(\mathbb{N}\cup\{\infty\})$. Since $\mathbb{N}\cup\{\infty\}$ is compact, $F_{R}(\mathbb{N}\cup\{\infty\})$ is topologically isomorphic to the free topological group $F_{M}(\mathbb{N}\cup\{\infty\})$. Therefore, $\pi_{1}^{top}(\Sigma((\mathbb{N}\cup\{\infty\})_{+}))$ is the free topological group on $\mathbb{N}\cup\{\infty\}$. The space $\Sigma(\mathbb{N}\cup\{\infty\})_{+}$ is illustrated here as a sequence of circles converging onto a "limit circle" all with a point in common.\begin{center} \includegraphics[height=4.5cm]{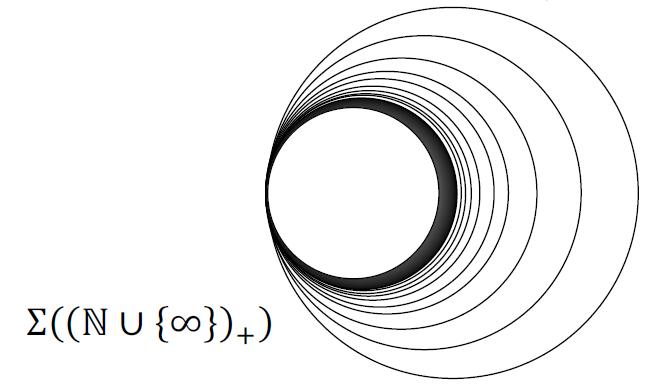} \end{center}
The topological fundamental group of this space gives a simple geometric interpretation of the topology of the free topological group $F_{M}(\mathbb{N}\cup\{\infty\})$. We can form similar geometric interpretations of free topological groups on other totally path disconnected, Tychonoff, $k_{\omega}$-spaces such as for the Cantor set.}
\end{example}
A rather surprising example of when $\pitop$ fails to be a topological group is the following
\begin{example} \label{rationalexample} \emph{
Let $\mathbb{Q}$ be the rational numbers with their usual topology so that the weak suspension space $w \Sigma (\mathbb{Q}_+)$ is a locally simply connected planar continuum. Since $\mathbb{Q}$ is totally path disconnected, $\pi_{1}^{top}(w \Sigma (\mathbb{Q}_+))\cong F_{R}(\mathbb{Q})$. In \cite{FOT} it is shown that $\mathbb{Q}$ fails to satisfy both conditions 3a. and 3b. of Theorem \ref{tychonoffconditions} and so $F_{R}(\mathbb{Q})$ is not a topological group. Therefore $\pi_{1}^{top}(w \Sigma (\mathbb{Q}_+))$ cannot be a topological group. }
\end{example}
Apparently the group $F_{R}(\mathbb{Q})$ fails to be a topological group due to the lack of some notion of "local compacness" in $\mathbb{Q}$. It seems then, that the topological fundamental group may fail to be a topological group due to some similar notion of "local compactness" in loop spaces. If a full solution to Question 3 becomes available, this informal description could perhaps be made more precise.
\subsection{Conclusions}
The computations and analysis of of sections 2 and 3 offer new insight into the nature of topological fundamental groups and provide a new geometric interpretation of many quasitopological and free topological groups. We note here how these ideas may be extended to higher dimensions and abelian groups, i.e. to the higher topological homotopy groups $\pi_{n}^{top}(X,x)=\pi_{0}^{top}(\Omega^{n}(X,x))$ and free abelian (Markov) topological groups. These groups were first mentioned in \cite{GHMM}, however, just as in \cite{Biss}, the authors falsely assert that products of quotient maps are quotient maps. The following question remains open.\\
\\
\textbf{Question 4:} For $n\geq 2$, is $\pi_{n}^{top}$ a functor to the category abelian topological groups?\\

As mentioned in the introduction, Fabel has shown that the topological fundamental group of the Hawaiian earring fails to be a topological group. In discussion, Fabel has indicated that the n-th topological fundamental group of the n-dimensional Hawaiian earring is indeed a topological group for each $n\geq 2$, this gives hope to an affirmative answer for Question 4. However, the results in this paper seem to indicate otherwise. It is well known (see \cite{Sipacheva}) that the abelian version of the reduction topology on free groups: $A_{R}(Y)$ is not always the free abelian topological group $A_M(Y)$ on $Y$ (the identity $A_{R}(Y)\ra A_{M}(Y)$ is continuous but not always open).\\
\\
\textbf{Conjecture:} $\pi_{n}^{top}(\Sigma^{n}(X_+))$ and $A_{R}(\piztop)$ are isomorphic as groups, however, the former has topology (sometimes strictly) finer than the later.\\

If this is indeed, the case $\pi_{n}^{top}$ will be a functor to the category of abelian quasitopological groups but not to the category of abelian topological groups. A computation of $\pi_{n}^{top}(\Sigma^{n}(X_+))$ for $n\geq 2$ should then provide an answer to Question 4.\\

Our application of free topological groups also indicates a potiential "fix" to the fundamental defect of $\pi_{1}^{top}$ which is shown to exist in this paper and in \cite{F3}: $\pitopx$ is not always a topological group even for some reasonable spaces $X$.\\

There are a number of alternative approaches to transfering topological structure to fundamental groups. Many of these approaches make use of some topology on the set of loops, however, if we change the topology on the set of loops $\Omega(X,x)$ and again identify homotopy classes, we will likely disturb the suspension-loop adjunction and be forced to give up many of the conveniences of classical homotopy. On the other hand, if we keep the compact-open topology, one might then ask if there is a finest group topology on $\pi_{1}(X,x)$ such that the canonical function $\Omega(X,x)\ra \pi_{1}(X,x)$ is continuous. Luckily, the existence of free topological groups gives such a topology for free: The function $\Omega(X,x)\ra \pi_{1}(X,x)$ induces a group epimorphism $F_{M}(\Omega(X,x))\ra \pi_{1}(X,x)$ and we may then view $\pi_{1}(X,x)$ as the quotient of the free topological group $F_{M}(\Omega(X,x))$. This construction will provide a homotopy invariant which is slightly weaker than $\pi_{1}^{top}$ but which takes values in the category of topological groups and has much nicer properties. A detailed study of this topology (at least when it is interesting) should amount to studying the properties of free topological groups on loop spaces of pathological spaces.
\section{Acknowledgements}
The author would like to thank his advisor Maria Basterra for her guidance and many helpful suggestions in writing this paper. Thanks are also due to the University of New Hampshire Graduate School for funding in the form of a summer TA and dissertation year fellowship.
\end{document}